\documentclass[12pt]{amsart}

\usepackage[matrix,arrow,curve,frame]{xy}
\usepackage{amsmath,amsthm,amssymb,enumerate}
\usepackage{latexsym}
\usepackage{amscd}
\usepackage[colorlinks=false]{hyperref}
\usepackage{euscript}

\newtheorem{thm}[subsection]{Theorem}
\newtheorem{defn}[subsection]{Definition}

\newtheorem{cor}[subsection]{Corollary}
\newtheorem{lemma}[subsection]{Lemma}

\newtheorem{remark}[subsection]{Remark}

\theoremstyle{definition}

\numberwithin{equation}{section}

\newfont{\german}{eufm10}

\begin{document}
\pagestyle{plain}

\title{Invariant theory and the $\mathcal{W}_{1+\infty}$ algebra with negative integral central charge}

\author{Andrew R. Linshaw}
\address{Fachbereich Mathematik, Technische Universit\"at Darmstadt, 64289 Darmstadt, Germany.}
\email{linshaw@mathematik.tu-darmstadt.de}
%\thanks{}

%\date{\today}
{\abstract

\noindent
The vertex algebra $\mathcal{W}_{1+\infty,c}$ with central charge $c$ may be defined as a module over the universal central extension of the Lie algebra of differential operators on the circle. For an integer $n\geq 1$, it was conjectured in the physics literature that $\mathcal{W}_{1+\infty,-n}$ should have a minimal strong generating set consisting of $n^2+2n$ elements. Using a free field realization of $\mathcal{W}_{1+\infty,-n}$ due to Kac-Radul, together with a deformed version of Weyl's first and second fundamental theorems of invariant theory for the standard representation of $GL_n$, we prove this conjecture. A consequence is that the irreducible, highest-weight representations of $\mathcal{W}_{1+\infty,-n}$ are parametrized by a closed subvariety of $\mathbb{C}^{n^2+2n}$.}

\maketitle
%\tableofcontents

\section{Introduction}

The Lie algebra $\mathcal{D}$ of regular differential operators on the circle has a universal central extension $\hat{\mathcal{D}} = \mathcal{D} \oplus \mathbb{C}\kappa$ which was introduced by Kac-Peterson in \cite{KP}. The representation theory of $\hat{\mathcal{D}}$ was first studied by Kac-Radul in \cite{KRI}, and in this paper the irreducible, quasi-finite highest-weight representations were constructed and classified. In \cite{FKRW}, the representation theory of $\hat{\mathcal{D}}$ was developed by Frenkel-Kac-Radul-Wang from the point of view of vertex algebras. For each $c\in\mathbb{C}$, $\hat{\mathcal{D}}$ admits a module $\mathcal{M}_c$ called the {\it vacuum module}, which is a vertex algebra freely generated by vertex operators $J^l(z)$, $l\geq 0$. The highest-weight representations of $\hat{\mathcal{D}}$ are in one-to-one correspondence with the highest-weight representations of $\mathcal{M}_c$. 

The unique irreducible quotient of $\mathcal{M}_c$ is a simple vertex algebra, and is often denoted by $\mathcal{W}_{1+\infty,c}$. These algebras have been studied extensively in the physics literature, and they also play an important role in the theory of integrable systems. Let $\pi_c$ denote the projection $\mathcal{M}_c\rightarrow \mathcal{W}_{1+\infty,c}$, whose kernel $\mathcal{I}_c$ is the maximal proper graded $\hat{\mathcal{D}}$-submodule of $\mathcal{M}_c$, and let $j^l = \pi_c(J^l)$. For $c\notin \mathbb{Z}$, $\mathcal{M}_c$ is irreducible, so $\mathcal{W}_{1+\infty,c}\cong \mathcal{M}_c$, but when $c$ is an integer $n$, $\mathcal{M}_{n}$ is reducible, and the structure and representation theory of $\mathcal{W}_{1+\infty,n}$ are nontrivial. 

For $n\geq 1$, $\mathcal{W}_{1+\infty,n}$ has a free field realization as the invariant space $\mathcal{E}(V)^{GL_n}$ \cite{FKRW}. Here $V=\mathbb{C}^n$, $\mathcal{E}(V)$ is the $bc$-system, or semi-infinite exterior algebra associated to $V$, and $\mathcal{E}(V)^{GL_n}$ is the invariant subalgebra under the natural action of $GL_n$ by vertex algebra automorphisms. Using this realization, the authors explicitly identified $\mathcal{W}_{1+\infty,n}$ with the vertex algebra $\mathcal{W}(\mathfrak{g}\mathfrak{l}_n)$ of central charge $n$, and classified its irreducible representations. In particular, $\mathcal{M}_n$ has a unique nontrivial singular vector (up to scalar multiples) of weight $n+1$, which generates $\mathcal{I}_n$ as a vertex algebra ideal. This singular vector gives rise to a \lq\lq decoupling relation" in $\mathcal{W}_{1+\infty,n}$ of the form $$j^n = P(j^0,\dots, j^{n-1}),$$ where $P$ is a normally ordered polynomial in the vertex operators $j^0,\dots,j^{n-1}$ and their derivatives. 

For $n\geq 1$, there is an analogous free field realization of $\mathcal{W}_{1+\infty,-n}$ as the invariant subalgebra $\mathcal{S}(V)^{GL_n}$, where $\mathcal{S}(V)$ is the $\beta\gamma$-system, or semi-infinite symmetric algebra, associated to $V=\mathbb{C}^n$ \cite{KRII}. In this paper, Kac-Radul used an infinite-dimensional version of the theory of Howe pairs to decompose $\mathcal{S}(V)$ into a direct sum of modules of the form $L\otimes M$ where, $L$ is an irreducible, finite-dimensional $GL_n$-module, and $M$ is an irreducible, highest-weight $\mathcal{W}_{1+\infty,-n}$-module. In particular, this decomposition of $\mathcal{S}(V)$ furnishes an interesting discrete set of irreducible, highest-weight $\mathcal{W}_{1+\infty,-n}$-modules. In \cite{A}, Adamovic used the realization $\mathcal{W}_{1+ \infty,-n}\cong \mathcal{S}(V)^{GL_n}$ together with the Friedan-Martinec-Shenker bosonization to exhibit $\mathcal{W}_{1+\infty,-n}$ as a subalgebra of the tensor product of $2n$ copies of the Heisenberg vertex algebra, and constructed a $2n$-parameter family of irreducible, highest-weight modules over $\mathcal{W}_{1+\infty,-n}$. However, in order to classify such modules, more information about the structure of $\mathcal{W}_{1+\infty,-n}$ and $\mathcal{I}_{-n}$ is needed, and these issues were not addressed in either of these papers. 

The first step in this direction was taken by Wang in \cite{WI}\cite{WII}. In the case $n=1$, he showed that $\mathcal{W}_{1+\infty,-1}$ is isomorphic to $\mathcal{W}(\mathfrak{g}\mathfrak{l}_3)$ with central charge $-2$, and classified its irreducible modules. He also conjectured in \cite{WIII} that $\mathcal{I}_{-1}$ should be generated by a unique singular vector.  However, for $n>1$, the structure of $\mathcal{W}_{1+\infty,-n}$ is still an open problem. There is a singular vector in $\mathcal{M}_{-n}$ of weight $(n+1)^2$, and it was conjectured in the physics literature by Blumenhagen-Eholzer-Honecker-Hornfeck-Hubel in \cite{B-H}, and also by Wang in \cite{WIII}, that this vector should give rise to a decoupling relation of the form \begin{equation} \label{dintro} j^l = P(j^0,\dots,j^{l-1}),\ \ \ \ \  l =n^2 +2n.\end{equation} The main result of this paper is a proof of this conjecture, and our starting point is the realization $\mathcal{W}_{1+\infty,-n}\cong \mathcal{S}(V)^{GL_n}$. This point of view allows us to study $\mathcal{W}_{1+\infty,-n}$ using {\it classical invariant theory}, an approach which was first suggested by Eholzer-Feher-Honecker in \cite{EFH}. As a vector space, $\mathcal{S}(V)^{GL_n}$, is isomorphic to the classical invariant ring $$R=(Sym \bigoplus_{k\geq 0} (V_k\oplus V^*_k))^{GL_n},$$ where $V_k$ and $V^*_k$ are copies of $V$ and $V^*$, respectively. We view $\mathcal{S}(V)^{GL_n}$ as a {\it deformation} of $R$, in the sense that $\mathcal{S}(V)^{GL_n}$ is linearly isomorphic to $R$, and admits a filtration for which the associated graded object $gr(\mathcal{S}(V)^{GL_n})$ is isomorphic to $R$ as a commutative ring. The generators and relations of $R$ are given by Weyl's first and second fundamental theorems of invariant theory for the standard representation of $GL_n$ \cite{W}. By a careful analysis of the deformation of this ring structure, we prove two key facts:
\begin{itemize}
\item For $n\geq 1$, $\mathcal{M}_{-n}$ has a unique nontrivial singular vector (up to scalar multiples) of weight $(n+1)^2$, which generates the maximal proper submodule $\mathcal{I}_{-n}$. This is analogous to the uniqueness of the singular vector in $\mathcal{M}_{n}$ of weight $n+1$ which generates $\mathcal{I}_n$, for $n\geq 1$.

\item This singular vector is of the form $J^l - P(J^0,\dots, J^{l-1})$ for $l=n^2+2n$, and hence gives rise to a decoupling relation in $\mathcal{W}_{1+\infty,-n}$ of the form (\ref{dintro}). Using this relation, it is easy to construct higher decoupling relations $j^r = Q_r(j^0,\dots, j^{l-1})$ for $r>l$. It follows that $\mathcal{W}_{1+\infty,-n}$ has a minimal strong generating set $\{j^0,\dots,j^{n^2+2n-1}\}$.

\end{itemize}

It is known \cite{FKRW} that the Zhu algebra of $\mathcal{M}_c$ is isomorphic to the polynomial algebra $\mathbb{C}[a^0,a^1,a^2,\dots]$. It follows from our main result that the Zhu algebra of $\mathcal{W}_{1+\infty,-n}$ is a quotient of the polynomial ring $\mathbb{C}[a^0,\dots,a^{n^2+2n-1}]$ by an ideal $I_{-n}$ corresponding to $\mathcal{I}_{-n}$. In particular, the Zhu algebra of $\mathcal{W}_{1+\infty,-n}$ is commutative, so its irreducible representations are one-dimensional, and are in one-to-one correspondence with the points on the variety $V(I_{-n})\subset \mathbb{C}^{n^2+2n}$. It follows that the irreducible, admissible representations of $\mathcal{W}_{1+\infty,-n}$ are all highest-weight representations, and are parametrized by $V(I_{-n})$ as well. We show that $V(I_{-n})$ is a proper, closed subvariety of $\mathbb{C}^{n^2+2n}$. In future work, we hope to study the geometry of this variety in more detail.

\section{Vertex algebras}
In this section, we define vertex algebras, which have been discussed from various different points of view in the literature \cite{B}\cite{FHL}\cite{FLM}\cite{K}\cite{LI}\cite{LZ}. We will follow the formalism developed in \cite{LZ} and partly in \cite{LI}. Let $V=V_0\oplus V_1$ be a super vector space over $\mathbb{C}$, and let $z,w$ be formal variables. By $QO(V)$, we mean the space of all linear maps $$V\rightarrow V((z)) =\{\sum_{n\in\mathbb{Z}} v(n) z^{-n-1}|
v(n)\in V,\ v(n)=0\ \text{for} \ n>>0 \}.$$ Each element $a\in QO(V)$ can be
uniquely represented as a power series
$$a=a(z) =\sum_{n\in\mathbb{Z}}a(n)z^{-n-1}\in (End\ V)[[z,z^{-1}]].$$ We
refer to $a(n)$ as the $n$th {\it Fourier mode} of $a(z)$. Each $a\in
QO(V)$ is assumed to be of the shape $a=a_0+a_1$ where $a_i:V_j\rightarrow V_{i+j}((z))$ for $i,j\in\mathbb{Z}/2\mathbb{Z}$, and we write $|a_i| = i$.

On $QO(V)$ there is a set of nonassociative bilinear operations
$\circ_n$, indexed by $n\in\mathbb{Z}$, which we call the $n$th {\it circle
products}. For homogeneous $a,b\in QO(V)$, they are defined by
$$
a(w)\circ_n b(w)=Res_z a(z)b(w)\ \iota_{|z|>|w|}(z-w)^n-
(-1)^{|a||b|}Res_z b(w)a(z)\ \iota_{|w|>|z|}(z-w)^n.
$$
Here $\iota_{|z|>|w|}f(z,w)\in\mathbb{C}[[z,z^{-1},w,w^{-1}]]$ denotes the
power series expansion of a rational function $f$ in the region
$|z|>|w|$. We usually omit the symbol $\iota_{|z|>|w|}$ and just
write $(z-w)^{-1}$ to mean the expansion in the region $|z|>|w|$,
and write $-(w-z)^{-1}$ to mean the expansion in $|w|>|z|$. It is
easy to check that $a(w)\circ_n b(w)$ above is a well-defined
element of $QO(V)$.

The nonnegative circle products are connected through the {\it
operator product expansion} (OPE) formula.
For $a,b\in QO(V)$, we have \begin{equation}\label{opeformula} a(z)b(w)=\sum_{n\geq 0}a(w)\circ_n
b(w)\ (z-w)^{-n-1}+:a(z)b(w):,\end{equation} which is often written as
$a(z)b(w)\sim\sum_{n\geq 0}a(w)\circ_n b(w)\ (z-w)^{-n-1}$, where
$\sim$ means equal modulo the term $$
:a(z)b(w):\ =a(z)_-b(w)\ +\ (-1)^{|a||b|} b(w)a(z)_+.$$ Here
$a(z)_-=\sum_{n<0}a(n)z^{-n-1}$ and $a(z)_+=\sum_{n\geq
0}a(n)z^{-n-1}$. Note that $:a(w)b(w):$ is a well-defined element of
$QO(V)$. It is called the {\it Wick product} of $a$ and $b$, and it
coincides with $a\circ_{-1}b$. The other negative circle products
are related to this by
$$ n!\ a(z)\circ_{-n-1}b(z)=\ :(\partial^n a(z))b(z): ,$$
where $\partial$ denotes the formal differentiation operator
$\frac{d}{dz}$. For $a_1(z),\dots ,a_k(z)\in QO(V)$, the $k$-fold
iterated Wick product is defined to be
\begin{equation}\label{iteratedwick} :a_1(z)a_2(z)\cdots a_k(z):\ =\ :a_1(z)b(z): ,\end{equation}
where $b(z)=\ :a_2(z)\cdots a_k(z): $. We often omit the formal variable $z$ when no confusion can arise.

The set $QO(V)$ is a nonassociative algebra with the operations
$\circ_n$ and a unit $1$. We have $1\circ_n a=\delta_{n,-1}a$ for
all $n$, and $a\circ_n 1=\delta_{n,-1}a$ for $n\geq -1$. A linear subspace $\mathcal{A}\subset QO(V)$ containing 1 which is closed under the circle products will be called a {\it quantum operator algebra} (QOA).
In particular $\mathcal{A}$ is closed under $\partial$
since $\partial a=a\circ_{-2}1$. Many formal algebraic
notions are immediately clear: a homomorphism is just a linear
map that sends $1$ to $1$ and preserves all circle products; a module over $\mathcal{A}$ is a
vector space $M$ equipped with a homomorphism $\mathcal{A}\rightarrow
QO(M)$, etc. A subset $S=\{a_i|\ i\in I\}$ of $\mathcal{A}$ is said to generate $\mathcal{A}$ if any element $a\in\mathcal{A}$ can be written as a linear
combination of nonassociative words in the letters $a_i$, $\circ_n$, for
$i\in I$ and $n\in\mathbb{Z}$. We say that $S$ {\it strongly generates} $\mathcal{A}$ if any $a\in\mathcal{A}$ can be written as a linear combination of words in the letters $a_i$, $\circ_n$ for $n<0$. Equivalently, $\mathcal{A}$ is spanned by the collection $\{ :\partial^{k_1} a_{i_1}(z)\cdots \partial^{k_m} a_{i_m}(z):|\ i_1,\dots,i_m\in I,\  k_1,\dots,k_m \geq 0\}$.

We say that $a,b\in QO(V)$ {\it quantum commute} if $(z-w)^N
[a(z),b(w)]=0$ for some $N\geq 0$. Here $[,]$ denotes the super bracket. This condition implies that $a\circ_n b = 0$ for $n\geq N$, so (\ref{opeformula}) becomes a finite sum. If $N$ can be chosen to be 0, we say that $a,b$ commute. A {\it commutative quantum operator algebra} (CQOA) is a QOA whose elements pairwise quantum commute. Finally, the notion of a CQOA is equivalent to the notion of a vertex algebra. Every CQOA $\mathcal{A}$ is itself a faithful $\mathcal{A}$-module, called the {\it left regular
module}. Define
$$\rho:\mathcal{A}\rightarrow QO(\mathcal{A}),\ \ \ \ a\mapsto\hat a,\ \ \ \ \hat
a(\zeta)b=\sum_{n\in\mathbb{Z}} (a\circ_n b)\ \zeta^{-n-1}.$$ Then $\rho$ is an injective QOA homomorphism,
and the quadruple of structures $(\mathcal{A},\rho,1,\partial)$ is a vertex
algebra in the sense of \cite{FLM}. Conversely, if $(V,Y,{\bf 1},D)$ is
a vertex algebra, the collection $Y(V)\subset QO(V)$ is a
CQOA. {\it We will refer to a CQOA simply as a
vertex algebra throughout the rest of this paper}.

The following are useful identities which measure the nonassociativity and noncommutativity of the Wick product, and the failure of the positive circle products to be derivations of the
Wick product. Let $a,b,c$ be vertex operators in some vertex algebra $\mathcal{A}$, and let $n > 0$. Then
\begin{equation}\label{vaidi} :(:ab:)c:-:abc:\ =\sum_{k\geq0}{1\over(k+1)!}\left(:(\partial^{k+1}a)(b\circ_k
c): +(-1)^{|a||b|}:(\partial^{k+1}b)(a\circ_k c):\right),\end{equation}
\begin{equation}\label{vaidii} :ab:-(-1)^{|a||b|}:ba:\ =\sum_{k\geq0}{(-1)^k\over(k+1)!}\partial^{k+1}(a\circ_kb),\end{equation}
\begin{equation}\label{vaidiii} a\circ_n(:bc:)-:(a\circ_nb)c:-(-1)^{|a||b|}:b(a\circ_nc):\ = \sum_{k=1}^n \binom{n}{k} (a\circ_{n-k}b)\circ_{k-1}c.\end{equation}

\section{Category $\mathcal{R}$}
In \cite{LL} we considered a certain category $\mathcal{R}$ of vertex algebras, together with a functor from $\mathcal{R}$ to the category of supercommutative rings. This functor provides a bridge between vertex algebras and commutative algebra, and it allows us to study vertex algebras $\mathcal{A}\in\mathcal{R}$ by using the tools of commutative algebra. 

\begin{defn} Let $\mathcal{R}$ be the category of vertex algebras $\mathcal{A}$ equipped with a $\mathbb{Z}_{\geq 0}$-filtration
\begin{equation}\label{goodi} \mathcal{A}_{(0)}\subset\mathcal{A}_{(1)}\subset\mathcal{A}_{(2)}\subset \cdots,\ \ \ \mathcal{A} = \bigcup_{k\geq 0}
\mathcal{A}_{(k)} \end{equation} such that $\mathcal{A}_{(0)} = \mathbb{C}$, and for all
$a\in \mathcal{A}_{(k)}$, $b\in\mathcal{A}_{(l)}$, we have
\begin{equation}\label{goodii} a\circ_n b\in\mathcal{A}_{(k+l)},\ \ \ \text{for}\
n<0,\end{equation}
\begin{equation}\label{goodiii} a\circ_n b\in\mathcal{A}_{(k+l-1)},\ \ \ \text{for}\
n\geq 0.\end{equation}
Elements $a(z)\in\mathcal{A}_{(d)}\setminus \mathcal{A}_{(d-1)}$ are said to have
degree $d$, and morphisms in $\mathcal{R}$ are vertex
algebra homomorphisms which preserve the filtration.\end{defn}

Filtrations on vertex algebras satisfying (\ref{goodii})-(\ref{goodiii}) were introduced in \cite{LII} and are known as {\it good increasing filtrations}. Setting $\mathcal{A}_{(-1)} = \{0\}$, the associated graded object $gr(\mathcal{A}) = \bigoplus_{k\geq 0}\mathcal{A}_{(k)}/\mathcal{A}_{(k-1)}$ is a $\mathbb{Z}_{\geq 0}$-graded associative, supercommutative algebra with a unit $1$ under a product induced by the Wick product on $\mathcal{A}$. In general, there is no natural linear map $\mathcal{A}\rightarrow gr (\mathcal{A})$, but for each $r\geq 1$ we have the projection \begin{equation}\label{proj} \phi_r: \mathcal{A}_{(r)} \rightarrow \mathcal{A}_{(r)}/\mathcal{A}_{(r-1)}\subset gr(\mathcal{A}).\end{equation}
Moreover, $gr(\mathcal{A})$ has a derivation $\partial$ of degree zero
(induced by the operator $\partial = \frac{d}{dz}$ on $\mathcal{A}$), and
for each $a\in\mathcal{A}_{(d)}$ and $n\geq 0$, the operator $a\circ_n$ on $\mathcal{A}$
induces a derivation of degree $d-k$ on $gr(\mathcal{A})$, which we also denote by $a\circ_n$. Here $$k  = sup \{ j\geq 1|\  \mathcal{A}_{(r)}\circ_n \mathcal{A}_{(s)}\subset \mathcal{A}_{(r+s-j)}\ \forall r,s,n\geq 0\},$$ as in \cite{LL}. Finally, these derivations give $gr(\mathcal{A})$ the structure of a vertex Poisson algebra.

The assignment $\mathcal{A}\mapsto gr(\mathcal{A})$ is a functor from $\mathcal{R}$ to the category of $\mathbb{Z}_{\geq 0}$-graded supercommutative rings with a differential $\partial$ of degree 0, which we will call $\partial$-rings. A $\partial$-ring is the same thing as an {\it abelian} vertex algebra, that is, a vertex algebra $\mathcal{V}$ in which $[a(z),b(w)] = 0$ for all $a,b\in\mathcal{V}$ \cite{B}. A $\partial$-ring $A$ is said to be generated by a subset $\{a_i|\ i\in I\}$ if $\{\partial^k a_i|\ i\in I, k\geq 0\}$ generates $A$ as a graded ring. The key feature of $\mathcal{R}$ is the following reconstruction property \cite{LL}:

\begin{lemma} \label{recon} Let $\mathcal{A}$ be a vertex algebra in $\mathcal{R}$ and let $\{a_i|\ i\in I\}$ be a set of generators for $gr(\mathcal{A})$ as a $\partial$-ring, where $a_i$ is homogeneous of degree $d_i$. If $a_i(z)\in\mathcal{A}_{(d_i)}$ are vertex operators such that $\phi_{d_i}(a_i(z)) = a_i$, then $\mathcal{A}$ is strongly generated as a vertex algebra by $\{a_i(z)|\ i\in I\}$.\end{lemma}

There is a similar reconstruction property for kernels of surjective morphisms in $\mathcal{R}$. Let $f:\mathcal{A}\rightarrow \mathcal{B}$ be a morphism in $\mathcal{R}$ with kernel $\mathcal{J}$, such that $f$ maps $\mathcal{A}_{(k)}$ onto $\mathcal{B}_{(k)}$ for all $k\geq 0$. The kernel $J$ of the induced map $gr(f): gr(\mathcal{A})\rightarrow gr(\mathcal{B})$ is a homogeneous $\partial$-ideal (i.e., $\partial J \subset J$). A set $\{a_i|\ i\in I\}$ such that $a_i$ is homogeneous of degree $d_i$ is said to generate $J$ as a $\partial$-ideal if $\{\partial^k a_i|\ i\in I,\ k\geq 0\}$ generates $J$ as an ideal.

\begin{lemma} \label{idealrecon} Let $\{a_i|\ i\in I\}$ be a generating set for $J$ as a $\partial$-ideal, where $a_i$ is homogeneous of degree $d_i$. Then there exist vertex operators $a_i(z)\in \mathcal{A}_{(d_i)}$ with $\phi_{d_i}(a_i(z)) = a_i$, such that $\{a_i(z)|\ i\in I\}$ generates $\mathcal{J}$ as a vertex algebra ideal.\end{lemma}
\begin{proof} First, let $a'_i(z)\in \mathcal{A}_{(d_i)}$ be an arbitrary vertex operator satisfying $\phi_{d_i}(a'_i(z)) = a_i$. Clearly $a'_i(z)$ need not lie in $\mathcal{J}$, but $f(a'_i(z))$ lies in $\mathcal{B}_{(d_i -1)}$. Since $f$ maps $\mathcal{A}_{(d_i-1)}$ onto $\mathcal{B}_{(d_i-1)}$, there exists $c_{i}(z)\in \mathcal{A}_{(d_i-1)}$ such that $f(c_{i}(z)) = -f(a'_i(z))$. Letting $a_i(z) = a'_i(z) + c_i(z)$, it follows that $a_i(z)\in \mathcal{J}$ and $\phi_{d_i}(a_i(z)) = a_i$. 

Now given $\omega(z)\in \mathcal{J}$ of degree $k$, we can write $\phi_k(\omega) = \sum_{i\in I}\sum_{j\geq 0} f_{ij} \partial^j a_i$, where all but finitely many $f_{ij}$ are zero, and each $f_{ij}$ is homogeneous of degree $k-d_i$. Choose vertex operators $f_{ij}(z)$ such that $\phi_{k-d_i}(f_{ij}(z)) = f_{ij}$, and let $$\omega'(z) = \sum_{i\in I}\sum_{j\geq 0} :f_{ij}(z) \partial^j a_i(z):.$$ Since each $a_i(z)\in \mathcal{J}$, $\omega''(z) = \omega(z) - \omega'(z)$ also lies in $\mathcal{J}$. Clearly $\phi_k(\omega'(z) - \omega(z)) = 0$, so $\omega''(z)\in \mathcal{J}\cap \mathcal{A}_{(k-1)}$. The claim follows by induction on $k$. \end{proof}

\section{The algebra $\mathcal{W}_{1+\infty,c}$}
Let $\mathcal{D}$ be the Lie algebra of regular differential operators on $\mathbb{C}\setminus \{0\}$, with coordinate $t$. A standard basis for $\mathcal{D}$ is $$J^l_k = -t^{l+k} (\partial_t)^l,\ \ \ k\in \mathbb{Z},\ \ \ l\in \mathbb{Z}_{\geq 0},$$ where $\partial_t = \frac{d}{dt}$. $\mathcal{D}$ has a 2-cocycle given by \begin{equation}\label{cocycle}\Psi\bigg(f(t) (\partial_t)^m,  g(t) (\partial_t)^n\bigg) = \frac{m! n!}{(m+n+1)!} Res_{t=0} f^{(n+1)}(t) g^{(m)}(t) dt,\end{equation} and a corresponding central extension $\hat{\mathcal{D}} = \mathcal{D} \oplus \mathbb{C} \kappa$, which was first studied by Kac-Peterson in \cite{KP}. $\hat{\mathcal{D}}$ has a $\mathbb{Z}$-grading $\hat{\mathcal{D}} = \bigoplus_{j\in\mathbb{Z}} \hat{\mathcal{D}}_j$ by weight, given by
$$wt (J^l_k) = k,\ \ \ \  wt (\kappa) = 0,$$ and a triangular decomposition $$\hat{\mathcal{D}} = \hat{\mathcal{D}}_+\oplus\hat{\mathcal{D}}_0\oplus \hat{\mathcal{D}}_-,$$ where $\hat{\mathcal{D}}_{\pm} = \bigoplus_{j\in \pm \mathbb{N}} \hat{\mathcal{D}}_j$ and $\hat{\mathcal{D}}_0 = \mathcal{D}_0\oplus \mathbb{C}\kappa.$
For a fixed $c\in\mathbb{C}$ and $\lambda\in \mathcal{D}_0^*$, define the Verma module with central charge $c$ over $\hat{\mathcal{D}}$ by
$$\mathcal{M}_c(\hat{\mathcal{D}},\lambda) = U(\hat{\mathcal{D}})\otimes_{U(\hat{\mathcal{D}}_0\oplus \hat{\mathcal{D}}_+)} C_{\lambda},$$ where $C_{\lambda}$ is the one-dimensional $\hat{\mathcal{D}}_0\oplus \hat{\mathcal{D}}_+$-module on which $\kappa$ acts by multiplication by $c$ and $h\in\hat{\mathcal{D}}_0$ acts by multiplication by $\lambda(h)$, and $\hat{\mathcal{D}}_+$ acts by zero. There is a unique irreducible quotient of $\mathcal{M}_c(\hat{\mathcal{D}},\lambda)$ denoted by $V_c(\hat{\mathcal{D}},\lambda)$.

Let $\mathcal{P}$ be the parabolic subalgebra of $\mathcal{D}$ consisting of differential operators which extend to all of $\mathbb{C}$, which has a basis $\{J^l_k|\ l\geq 0,\ l+k\geq 0\}$. The cocycle $\Psi$ vanishes on $\mathcal{P}$, so $\mathcal{P}$ may be regarded as a subalgebra of $\hat{\mathcal{D}}$. Clearly $\hat{\mathcal{D}}_0\oplus \hat{\mathcal{D}}_+\subset \hat{\mathcal{P}}$, where $\hat{\mathcal{P}} = \mathcal{P}\oplus \mathbb{C}\kappa$. The induced $\hat{\mathcal{D}}$-module $$\mathcal{M}_c=\mathcal{M}_c(\hat{\mathcal{D}},\hat{\mathcal{P}}) = U(\hat{\mathcal{D}})\otimes_{U(\hat{\mathcal{P}})} C_0$$ is then a quotient of $\mathcal{M}_c(\hat{\mathcal{D}},0)$, and is known as the {\it vacuum $\hat{\mathcal{D}}$-module of central charge $c$}. $\mathcal{M}_c$ has the structure of a vertex algebra which is generated by fields $$J^l(z) = \sum_{k\in\mathbb{Z}} J^l_k z^{-k-l-1},\ \ \ l\geq 0$$ of weight $l+1$. The modes $J^l_k$ represent $\hat{\mathcal{D}}$ on $\mathcal{M}_c$, and in order to be consistent with our earlier notation, we rewrite these fields in the form
$$J^l(z) = \sum_{k\in\mathbb{Z}} J^l(k) z^{-k-1},$$ where $J^l(k) = J^l_{k-l}$. In fact, $\mathcal{M}_c$ is {\it freely} generated by $\{J^l(z)|\ l\geq 0\}$; the set of iterated Wick products \begin{equation}\label{standmon} :\partial^{i_1}J^{l_1}(z)\cdots \partial^{i_r} J^{l_r}(z):,\end{equation} such that $l_1\leq \cdots \leq l_r$ and $i_a\leq i_b$ if $l_a = l_b$, forms a basis for $\mathcal{M}_c$. Define a filtration $$(\mathcal{M}_c)_{(0)} \subset (\mathcal{M}_c)_{(1)}\subset \cdots$$ on $\mathcal{M}_c$ as follows: for $k\geq 0$, $(\mathcal{M}_c)_{(2k)}$ is the span of monomials of the form (\ref{standmon}), for $r\leq k$, and $(\mathcal{M}_c)_{(2k+1)} = (\mathcal{M}_c)_{(2k)}$. In particular, each $J^l$ and its derivatives has degree $2$. Equipped with this filtration, $\mathcal{M}_c$ lies in the category $\mathcal{R}$, and $gr(\mathcal{M}_c)$ is the polynomial algebra $\mathbb{C}[\partial^k J^l|\  k,l\geq 0]$. Moreover, the vertex Poisson algebra structure on $\mathcal{M}_c$ is the same for all $c$. In particular, each operator $J^l\circ_k$ for $k,l\geq 0$ is a derivation of degree zero on $gr(\mathcal{M}_c)$, and this action of $\mathcal{P}$ on $gr(\mathcal{M}_c)$ is independent of $c$.

\begin{lemma} \label{weakfg} For each $c\in\mathbb{C}$, $\mathcal{M}_c$ is generated as a vertex algebra by  $J^0$, $J^1$, and $J^2$.\end{lemma} 

\begin{proof} Let $\mathcal{J}$ be the vertex subalgebra of $\mathcal{M}_c$ generated by $J^0$, $J^1$, and $J^2$. An OPE calculation shows that for $l\geq 1$, $$J^2 \circ_1 J^{l-1} = -(l+1) J^{l} +2 \partial J^{l-1},\ \ \ \ \ \ \ J^1\circ_0 J^l= -\partial J^{l}.$$
It follows that $\alpha\circ_1 J^{l-1} = -(l+1) J^l$, where $\alpha = J^2 -2 \partial J^1$.  Since $\alpha\in \mathcal{J}$, it follows by induction that $J^l\in \mathcal{J}$ for all $l$. \end{proof}

In particular, $\mathcal{M}_c$ is a finitely generated vertex algebra. However, $\mathcal{M}_c$ is not {\it strongly} generated by any finite set of vertex operators. This follows from the fact that $gr(\mathcal{M}_c)\cong \mathbb{C}[\partial^k J^l|\ k,l\geq 0]$, which implies that there are no normally ordered polynomial relations among the vertex operators $J^l$, $l \geq 0$, and their derivatives.

A weight-homogeneous element $\omega\in \mathcal{M}_c$ is called a {\it singular vector} if $J^l\circ_k \omega = 0$ for all $k>l\geq 0$. The maximal proper $\hat{\mathcal{D}}$-submodule $\mathcal{I}_c$ is the vertex algebra ideal generated by all singular vectors $\omega\neq 1$, and the unique irreducible quotient $\mathcal{M}_c/\mathcal{I}_c$ is often denoted by $\mathcal{W}_{1+\infty,c}$ in the literature. We denote the projection $\mathcal{M}_c\rightarrow \mathcal{W}_{1+\infty,c}$ by $\pi_{c}$, and we use the notation \begin{equation}\label{lowercasej} j^l = \pi_c (J^l),\ \ \ \ \ \ l\geq 0\end{equation} in order to distinguish between $J^l\in \mathcal{M}_{c}$ and its image in $\mathcal{W}_{1+\infty,c}$ Clearly $\mathcal{W}_{1+\infty,c}$ is generated by $j^0,j^1,j^2$ as a vertex algebra, but there may now be normally ordered polynomial relations among $\{\partial^k j^l|\ k,l\geq 0\}$.

For $c\notin \mathbb{Z}$, $\mathcal{M}_c$ is irreducible, so $\mathcal{W}_{1+\infty,c} = \mathcal{M}_c$, but for $n\in\mathbb{Z}$, $\mathcal{M}_n$ is reducible. For $n \geq 1$, $\mathcal{W}_{1+\infty,n}$ is known to be isomorphic to $\mathcal{W}(\mathfrak{g}\mathfrak{l}_n)$ with central charge $n$ \cite{FKRW}. An important ingredient in the proof is a realization of $\mathcal{W}_{1+\infty,n}$ as a subalgebra of the $bc$-system $\mathcal{E}(V)$, for $V=\mathbb{C}^n$. This vertex algebra was introduced in \cite{FMS}, and is the unique vertex algebra with odd generators $b^x(z), c^{x'}(z)$ for $x\in V$ and $x'\in V^*$, satisfying the OPE relations $$b^x(z) c^{x'}(w)\sim\langle x',x\rangle (z-w)^{-1},\ \ \ \ \ \ c^{x'}(z)b^x(w)\sim \langle x',x\rangle (z-w)^{-1},$$ \begin{equation}\label{bcope} b^x(z)b^y(w)\sim 0,\ \ \ \ \ \ c^{x'}(z)c^{y'}(w)\sim 0,\end{equation} where $\langle,\rangle$ denotes the natural pairing between $V^*$ and $V$. The map \begin{equation}\label{bcrealize} \mathcal{W}_{1+\infty,n}\rightarrow \mathcal{E}(V),\ \ \ \ \ \ j^l(z)\mapsto \sum_{i=1}^n : c^{x'_i}(z) \partial^l b^{x_i}(z):\end{equation} identifies $\mathcal{W}_{1+\infty,n}$ with the invariant subalgebra $\mathcal{E}(V)^{GL_n}$ \cite{FKRW}. Here $\{x_1,\dots,x_n\}$ is a basis for $V$ and $\{x'_1,\dots,x'_n\}$ is the dual basis for $V^*$. There is a singular vector in $\mathcal{M}_n$ of weight $n+1$, which generates $\mathcal{I}_n$ and gives rise to a decoupling relation of the form $$j^n = P(j^0,\dots,j^{n-1}),$$ where $P$ is a normally ordered polynomial in the vertex operators $j^0,\dots,j^{n-1}$ and their derivatives. From this relation, higher decoupling relations $j^r = Q_r(j^0,\dots, j^{n-1})$ can be constructed for all $r>n$, so $\mathcal{W}_{1+\infty,n}$ is strongly generated by $\{j^0,\dots, j^{n-1}\}$.

For $n=-1$, $\mathcal{W}_{1+\infty,-1}$ is isomorphic to $\mathcal{W}(\mathfrak{g}\mathfrak{l}_3)$ with central charge $-2$, which was shown using the Friedan-Martinec-Shenker bosonization \cite{WI}\cite{WII}. In particular, $W_{1+\infty,-1}$ is a tensor product of a Heisenberg algebra $\mathcal{H}$ and the simple Zamolodchikov $\mathcal{W}_3$ algebra with $c=-2$. However, for $n>1$, the structure of $\mathcal{W}_{1+\infty,-n}$ is not well understood. There is a singular vector in $\mathcal{M}_{-n}$ of weight $(n+1)^2$, and it was conjectured in \cite{B-H} and \cite{WIII} that this singular vector gives rise to a decoupling relation \begin{equation}\label{decoupling} j^{l}= P(j^0,\dots, j^{l-1}) \end{equation} for $l=n^2+2n$. From such a relation, one can construct higher decoupling relations $j^r = Q_r(j^0,\dots, j^{l-1})$ for all $r>l$, so this would imply that $\mathcal{W}_{1+\infty,-n}$ is strongly generated by $\{j^0,\dots,j^{n^2+2n-1}\}$.

For $n\geq 1$, $\mathcal{W}_{1+\infty,-n}$ has an analogous realization as a subalgebra of the $\beta\gamma$-system $\mathcal{S}(V)$ for $V=\mathbb{C}^n$. The $\beta\gamma$-system, or algebra of chiral differential operators on $V$, was introduced in \cite{FMS}. It is the unique vertex algebra with even generators $\beta^{x}(z)$, $\gamma^{x'}(z)$ for $x\in V$, $x'\in V^*$, which satisfy the OPE relations $$\beta^x(z)\gamma^{x'}(w)\sim\langle x',x\rangle (z-w)^{-1},\ \ \ \ \ \ \gamma^{x'}(z)\beta^x(w)\sim -\langle x',x\rangle (z-w)^{-1},$$ \begin{equation}\label{betagamma} \beta^x(z)\beta^y(w)\sim 0,\ \ \ \ \ \ \gamma^{x'}(z)\gamma^{y'}(w)\sim 0.\end{equation} We give $\mathcal{S}(V)$ the conformal structure \begin{equation}\label{virasoroelement} L(z) = \sum_{i=1}^n :\beta^{x_i}(z)\partial\gamma^{x'_i}(z):,\end{equation} under which $\beta^{x_i}(z)$ and $\gamma^{x'_i}(z)$ are primary of conformal weights $1$ and $0$, respectively. Moreover, $\mathcal{S}(V)$ has a basis consisting of the normally ordered monomials
\begin{equation}\label{basisofs} :\partial^{I_1} \beta^{x_1}\cdots \partial^{I_n} \beta^{x_n}\partial^{J_1} \gamma^{x'_1}\cdots\partial^{J_n} \gamma^{x'_n}: .\end{equation} In this notation, $I_k = (i^k_1,\dots, i^k_{r_k})$ and $J_k =(j^k_1,\dots, j^k_{s_k})$ are lists of integers satisfying $0\leq i^k_1\leq \cdots \leq i^k_{r_k}$ and  $0\leq j^k_1\leq \cdots \leq j^k_{s_k}$, and 
$$\partial^{I_k}\beta^{x_k} = \ :\partial^{i^k_1}\beta^{x_k} \cdots \partial^{i^k_{r_k}} \beta^{x_k}:,\ \ \ \ \ \ \partial^{J_k}\gamma^{x'_k} = \ :\partial^{j^k_1}\gamma^{x'_k} \cdots \partial^{j^k_{s_k}} \gamma^{x'_k}:.$$ $\mathcal{S}(V)$ then has a $\mathbb{Z}_{\geq 0}$-grading \begin{equation}\label{grading} \mathcal{S}(V) = \bigoplus_{d\geq 0} \mathcal{S}(V)^{(d)},\end{equation} where $\mathcal{S}(V)^{(d)}$ is spanned by monomials of the form (\ref{basisofs}) of total degree $ d = \sum_{k=1}^n r_k + s_k$. Finally, we define the filtration $\mathcal{S}(V)_{(d)} = \bigoplus_{i=0}^d \mathcal{S}(V)^{(i)}$. This filtration satisfies (\ref{goodi})-(\ref{goodiii}), and we have
$$gr(\mathcal{S}(V))\cong Sym(\bigoplus_{k\geq 0} (V_k \oplus V^*_k)\big),\ \ \ \ \ V_k = \{\beta^{x}_k |\  x\in V\},\ \ \ \ \ V^*_k = \{\gamma^{x'}_k |\  x'\in V^*\}.$$ In this notation, $\beta^{x}_k$ and $\gamma^{x'}_k$ are the images of $\partial^k \beta^{x}(z)$ and $\partial^k\gamma^{x'}(z)$ in $gr(\mathcal{S}(V))$ under the projection $\phi_1: \mathcal{S}(V)_{(1)}\rightarrow \mathcal{S}(V)_{(1)}/\mathcal{S}(V)_{(0)}$. The embedding $\mathcal{W}_{1+\infty,-n}\rightarrow \mathcal{S}(V)$ introduced in \cite{KRII} is defined by
\begin{equation}\label{bgrealization} j^l(z) \mapsto \sum_{i=1}^n :\gamma^{x'_i}(z) \partial^l \beta^{x_i}(z):.\end{equation}
This map preserves conformal weight, and is a morphism in the category $\mathcal{R}$. For the rest of this paper, we will identify $\mathcal{W}_{1+\infty,-n}$ with its image in $\mathcal{S}(V)$.

The standard representation of $GL_n$ on $V = \mathbb{C}^n$ induces an action of $GL_n$ on $\mathcal{S}(V)$ by vertex algebra automorphisms. In fact, $GL_n$ is the {\it full} automorphism group of $\mathcal{S}(V)$ preserving the conformal structure (\ref{virasoroelement}). As shown in \cite{KRII}, $\mathcal{W}_{1+\infty,-n}$ is precisely the invariant subalgebra $\mathcal{S}(V)^{GL_n}$. The action of $GL_n$ on $\mathcal{S}(V)$ preserves the grading (\ref{grading}), so that $\mathcal{W}_{1+\infty,-n}$ is a graded subalgebra of $\mathcal{S}(V)$. We write \begin{equation}\label{gradingw} \mathcal{W}_{1+\infty,-n} = \bigoplus_{d\geq 0} ( \mathcal{W}_{1+\infty,-n})^{(d)},\ \ \ \ \ \ (\mathcal{W}_{1+\infty,-n})^{(d)} = \mathcal{W}_{1+\infty,-n}\cap\mathcal{S}(V)^{(d)},\end{equation} and define the corresponding filtration by $(\mathcal{W}_{1+\infty,-n} )_{(d)} = \bigoplus_{i=0}^{d} (\mathcal{W}_{1+\infty,-n} )^{(i)}$. 

The identification $\mathcal{W}_{1+\infty,-n}\cong \mathcal{S}(V)^{GL_n}$ suggests an alternative strong generating set for $\mathcal{W}_{1+\infty,-n}$. Define \begin{equation}\label{newgen} \omega_{a,b}(z) = \sum_{i=1}^n :\partial^a\beta^{x_i}(z) \partial ^b \gamma^{x'_i}(z):,\ \ \ \ \ a,b\geq 0.\end{equation} For each $m\geq 0$, let $A_m$ denote the vector space of dimension $m+1$ with basis $\{\omega_{a,b}|\  a+b = m\}$. Note that $\partial \omega_{a,b} = \omega_{a+1,b} + \omega_{a,b+1}$, so $\partial(A_m)\subset A_{m+1}$. Moreover, $A_m / \partial(A_{m-1})$ is one-dimensional, and $j^m = \omega_{m,0}\notin \partial (A_{m-1})$, so we have a decomposition \begin{equation}\label{decompofa} A_m = \partial A_{m-1}\oplus \langle j^m\rangle ,\end{equation} where $\langle j^m\rangle$ is the linear span of $j^m$. Finally, $\{\partial^a j^b|\  a+b = m\}$ is another basis for $A_m$, so each $\omega_{a,b}\in A_m$ can be expressed uniquely in the form
\begin{equation}\label{lincomb} \omega_{a,b} =\sum_{i=0}^m c_i \partial^i j^{m-i}\end{equation} for constants $c_0,\dots, c_m$. Hence $\{\omega_{a,b}|\ a,b\geq 0\}$ is another strong generating set for $\mathcal{W}_{1+\infty,-n}$ as a vertex algebra. The formula (\ref{lincomb}) holds in $\mathcal{W}_{1+\infty,-n}$ for any $n$, and allows us to define a new generating set $\{\Omega_{a,b}|\ a,b\geq 0\}$ for $\mathcal{M}_{-n}$, where $$\Omega_{a,b} =\sum_{i=0}^m c_i \partial^i J^{m-i},\ \ \ \ \ \ a+b = m.$$ In fact, this new generating set makes sense in $\mathcal{M}_c$ for any $c\in\mathbb{C}$, and $\pi_c(\Omega_{a,b}) = \omega_{a,b}$. Since $gr(\mathcal{M}_c) \cong \mathbb{C}[\partial^k J^l|\ k,l\geq 0]$, we also have $gr(\mathcal{M}_c) \cong \mathbb{C}[\Omega_{a,b}|\ a,b\geq 0]$. We will use the same notation $A_m$ to denote the linear span of $\{\Omega_{a,b}|\  a+b = m\}$, when no confusion can arise.

As shown in \cite{KRII}, the generating set (\ref{newgen}) has a natural interpretation in terms of Weyl's description of the ring of polynomial invariants for the standard representation of $GL_n$ \cite{W}.

\begin{thm}(Weyl) \label{weylfft} For $k\geq 0$, let $V_k$ be the copy of the standard $GL_n$-module $\mathbb{C}^n$ with basis $x_{i,k}$ for $i=1,\dots,n$, and let $V^*_k$ be the copy of $V^*$ with basis $x'_{i,k}$, $i=1,\dots,n$. The invariant ring $(Sym \bigoplus_{k\geq 0} (V_k\oplus V^*_k))^{GL_n}$ is generated by the quadratics \begin{equation}\label{weylgenerators} q_{a,b} = \sum_{i=1}^n x_{i,a} x'_{i,b},\end{equation} which correspond to the $GL_n$-invariant pairings $V_a\otimes V^*_b\rightarrow \mathbb{C}$ for $a,b\geq 0$. Let $Q_{a,b}$ be commuting indeterminates for $a,b\geq 0$. The kernel $I_n$ of the homomorphism \begin{equation}\label{weylquot} \mathbb{C}[Q_{a,b}]\rightarrow (Sym \bigoplus_{k\geq 0} (V_k\oplus V^*_k))^{GL_n},\ \ \ \ \ \ Q_{a,b}\mapsto q_{a,b},\end{equation} is generated by the $(n+1)\times (n+1)$ determinants \begin{equation}\label{weylrel} d_{I,J} = \det \left[\begin{matrix} Q_{i_0,j_0} & \cdots & Q_{i_0,j_n} \cr  \vdots  & & \vdots  \cr  Q_{i_n,j_0}  & \cdots & Q_{i_n,j_n} \end{matrix}\right].\end{equation} Here $I=(i_0,\dots, i_{n})$ and $J = (j_0,\dots, j_{n})$ are lists of integers satisfying \begin{equation}\label{ijineq} 0\leq i_0<\cdots <i_n,\ \ \ \ \ \  0\leq j_0<\cdots <j_n.\end{equation} \end{thm}

Since the action of $GL_n$ on $\mathcal{S}(V)$ preserves the filtration, we have $gr(\mathcal{S}(V)^{GL_n}) \cong (gr(\mathcal{S}(V))^{GL_n} \cong (Sym\bigoplus_{k\geq 0} (V_k\oplus V^*_k))^{GL_n}$, and under the projection $$\phi_2: \mathcal{S}(V)_{(2)}\rightarrow \mathcal{S}(V)_{(2)}/\mathcal{S}(V)_{(1)}\subset gr(\mathcal{S}(V)),$$ $\omega_{a,b}$ corresponds to $q_{a,b}$. 

Recall that the projection $\pi_{-n}: \mathcal{M}_{-n}\rightarrow \mathcal{W}_{1+\infty,-n}$ sending $\Omega_{a,b}\mapsto \omega_{a,b}$ is a morphism in the category $\mathcal{R}$. Under the identifications $$gr(\mathcal{M}_{-n})\cong \mathbb{C}[Q_{a,b}|\ a,b\geq 0],\ \ \ \ gr(\mathcal{W}_{1+\infty,-n})\cong Sym(\bigoplus (V_k\oplus V^*_k))^{GL_n}\cong \mathbb{C}[q_{a,b}]/I_n,$$ $gr(\pi_{-n})$ is just the quotient map (\ref{weylquot}). Clearly $\pi_{-n}$ maps each filtered piece $(\mathcal{M}_{-n})_{(k)}$ onto $(\mathcal{W}_{1+\infty,-n})_{(k)}$, so the hypotheses of Lemma \ref{idealrecon} are satisfied. Since $I_{n} = Ker (gr(\pi_{-n}))$ is generated by the determinants $d_{I,J}$, we can apply Lemma \ref{idealrecon} to find vertex operators $D_{I,J}\in (\mathcal{M}_{-n})_{(2n+2)}$ satisfying $\phi_{2n+2}(D_{I,J}) = d_{I,J}$, such that $\{D_{I,J}\}$ generates $\mathcal{I}_{-n}$. Since $\Omega_{a,b}$ has weight $a+b+1$, it follows that \begin{equation}\label{wtod} wt(D_{I,J}) = |I| + |J| +n+1,\ \ \ \ \ \ |I| =\sum_{a=0}^n i_a,\ \ \ \ \ \ |J| = \sum_{a=0}^n j_a.\end{equation}
In general, the vertex operators $a_i(z)$ furnished by Lemma \ref{idealrecon} satisfying $\phi_{d_i}(a_i(z)) = a_i$ which generate $\mathcal{I}$ are not unique. However, in our case, $D_{I,J}$ is uniquely determined by the conditions \begin{equation}\label{uniquedij} \phi_{2n+2}(D_{I,J}) = d_{I,J},\ \ \ \ \ \pi_{-n}(D_{I,J}) = 0.\end{equation} To see this, suppose that $D'_{I,J}$ is another vertex operator satisfying (\ref{uniquedij}). Then $D_{I,J} - D'_{I,J}$ lies in $(\mathcal{M}_{-n})_{(2n)} \cap \mathcal{I}_{-n}$, and since there are no relations in $\mathcal{W}_{1+\infty,-n}$ of degree less than $2n+2$, we have $D_{I,J} - D'_{I,J}=0$.

For each $n\geq 1$, define \begin{equation}\label{vsu} U_n = (\mathcal{M}_{-n})_{(2n+2)}\cap \mathcal{I}_{-n},\end{equation} which is just the vector space with basis $\{D_{I,J}\}$, where $I,J$ satisfy (\ref{ijineq}).

\begin{lemma} For all $n\geq 1$, $U_n$ is a module over the parabolic Lie algebra $\mathcal{P}\subset \hat{\mathcal{D}}$ generated by $\{J^l(k) = J^l\circ_k |\ k,l\geq 0\}$. \end{lemma}

\begin{proof} The action of $\mathcal{P}$ preserves the filtration degree, and in particular preserves $(\mathcal{M}_{-n})_{(2n+2)}$. Also, $\mathcal{P}$ preserves $\mathcal{I}_{-n}$ since $\mathcal{I}_{-n}$ is a vertex algebra ideal. \end{proof}

It will be convenient to work in the basis 
$\{\Omega_{a,b}\circ_{a+b-w}|\  a,b\geq 0,\  a+b-w\geq 0\}$ for $\mathcal{P}$. Note that $\Omega_{a,b}\circ_{a+b-w}$ is homogeneous of weight $w$. The action of $\mathcal{P}$ by derivations of degree zero on $gr(\mathcal{M}_{-n})$ coming from the vertex Poisson algebra structure is independent of $n$, and is specified by the action of $\mathcal{P}$ on the generators $\Omega_{l,m}$. We compute \begin{equation}\label{actionp} \Omega_{a,b}\circ_{a+b-w} \Omega_{l,m} = \lambda_{a,b,w,l} (\Omega_{l+w,m}) + \mu_{a,b,w,m} ( \Omega_{l,m+w}),\end{equation} where $$ \lambda_{a,b,w,l}  =  \bigg\{ \begin{matrix} (-1)^{b+1} \frac{(b+l)!}{(l+w-a)!} & l+w-a \geq 0 \cr & \cr 0 & l+w-a <0 \end{matrix}\ ,\ \ \ \   \mu_{a,b,w,m}  =  \bigg\{ \begin{matrix} (-1)^{a} \frac{(a+m)!}{(m+w-b)!} & m+w-b \geq 0 \cr & \cr 0 & m+w-b <0 \end{matrix}.$$

The action of $\mathcal{P}$ on $U_n$ is by \lq\lq weighted derivation" in the following sense. Fix $I = (i_0,\dots,i_n)$ and $J = (j_0,\dots,j_n)$, and let $D_{I,J}\in U_n$ be the corresponding element. Given $p = \Omega_{a,b}\circ_{a+b-w}\in \mathcal{P}$, we have \begin{equation}\label{paraction} p(D_{I,J}) = \sum_{r=0}^n \lambda_r D_{I^r,J} + \sum_{r=0}^n \mu_r D_{I,J^r},\end{equation} for lists $I^r = (i_0,\dots, i_{r-1}, i_r + w,i_{r+1},\dots, i_n)$ and $J^r = (j_0,\dots, j_{r-1}, j_r+ w,j_{r+1},\dots, j_n)$, and constants $\lambda_r$, $\mu_r$. If $i_r + w$ appears elsewhere on the list $I^r$, $\lambda_r = 0$, and if $j_r + w$ appears elsewhere on the list $J^r$, $\mu_r = 0$. Otherwise, \begin{equation}\label{actioni} \lambda_r = \pm \lambda_{a,b,w,i_r},\ \ \ \ \ \ \mu_r = \pm \mu_{a,b,w,j_r},\end{equation} where the signs $\pm$ are the signs of the permutations transforming $I^r$ and $J^r$ into lists in increasing order, as in (\ref{ijineq}).

For each $n\geq 1$, there is a distinguished element $D_0\in U_n$, defined by \begin{equation}\label{defodo} D_0 = D_{I,J},\ \ \ \ \ \ I = (0,1,\dots,n),\ \ \ \ \ \ J = (0,1,\dots,n).\end{equation} In this case, $|I| = |J| = \frac{n(n+1)}{2}$, so $D_0$ is the unique element of $\mathcal{I}_{-n}$ of minimal weight $(n+1)^2$, and hence is a singular vector in $\mathcal{M}_{-n}$. 

\begin{thm} \label{uniquesv} The element $D_0$ generates the ideal $\mathcal{I}_{-n}$. It follows that $D_0$ is the unique nontrivial singular vector in $\mathcal{M}_{-n}$, up to scalar multiples.\end{thm} 

\begin{proof} Since $\mathcal{I}_{-n}$ is generated by $U_n$ as a vertex algebra ideal, it suffices to show that $U_n$ is generated by $D_0$ as a module over $\mathcal{P}$. Let $U_n[k]$ denote the subspace of $U_n$ of weight $k$. Note that $U_n[k]$ is trivial for $k<(n+1)^2$ and is spanned by $D_0$ for $k = (n+1)^2$.

For $k>(n+1)^2$, we define a property $P_{n,k}$ of the subspace $U_n[k]$ as follows:

\begin{itemize}
\item For every $D_{I,J} \in U_n[k]$, if $I \neq (0,\dots,n)$, there exists an integer $s>0$ and elements $$D_{I^1,J},\dots, D_{I^s,J} \in U_n[k-1],$$ with $I^r = (i^r_0,\dots, i^r_n)$ satisfying $0\leq i_0 - i^r_0 \leq 1 $, and there exist elements $p_1,\dots,p_s\in\mathcal{P}$ satisfying
\begin{equation}\label{eltpi} \sum_{r=1}^s  p_r (D_{I^r,J})  = D_{I,J}.\end{equation}

\item Each $p_r$ appearing in (\ref{eltpi}) is a linear combination of elements of the form $\Omega_{a,b}\circ_{a+b-1}$, where $a\geq i_0$, and $b$ can be arbitrarily large. In other words, for each integer $t\geq 0$, we can assume that $p_1,\dots,p_s$ are linear combinations of elements $\Omega_{a,b}\circ_{a+b-1}$ with $b>t$. 

\item Similarly, if $I = (0,\dots,n)$ and $J = (j_0,\dots, j_n)\neq (0,\dots,n)$, there exists an integer $s>0$ and elements $D_{I,J^1},\dots, D_{I,J^s} \in U_n[k-1]$, with $J^r = (j^r_0,\dots, j^r_n)$, satisfying $0\leq j_0 - j^r_0\leq 1$, and elements $p_1,\dots,p_s\in\mathcal{P}$ satisfying
\begin{equation}\label{eltpii} \sum_{r=1}^s p_r (D_{I,J^r})= D_{I,J}.\end{equation} Moreover, each $p_r$ in (\ref{eltpii}) is a linear combination of elements of the form $\Omega_{a,b}\circ_{a+b-1}$, where $b\geq j_0$, and $a$ can be arbitrarily large. \end{itemize}

In order to prove that $U_n$ is generated by $D_0$ as a $\mathcal{P}$-module, it suffices to prove that property $P_{n,k}$ holds for all $n\geq 1$ and all $k>(n+1)^2$. We proceed by induction on $n$ and $k$. First, we show that $P_{n,k}$ holds for arbitrary $n$, and $k = (n+1)^2 + 1$. Second, we fix $n=1$ and show that property $P_{1,k}$ holds for all $k>4$ by induction on $k$. Third, we assume inductively that property $P_{m,k}$ holds for each pair $(m,k)$ with $m<n$ and $k > (m+1)^2$, and show that $P_{n,k}$ must then hold for $k>(n+1)^2$.

\smallskip

{\bf Step 1}: Fix $n\geq 1$ and let $k = (n+1)^2 + 1$. Then $U_n[k]$ is spanned by two elements corresponding to $\{I,J\} = \{(0,\dots,n-1,n+1),(0,\dots,n)\}$ and $\{I,J\} = \{(0,\dots,n),(0,\dots,n-1,n+1)\}$.

{\it Case 1}: $\{I,J\} = \{(0,\dots,n-1,n+1),(0,\dots,n)\}$. Let $I' = (0,\dots,n)$, so that $D_{I',J} = D_0$. For any $b>n+1$, let $$p = (-1)^{b+1}\frac{(n+1)!}{(b+n)!} \ \Omega_{0,b}\circ_{b-1}.$$ It follows from (\ref{actionp}), (\ref{paraction}), and (\ref{actioni}) that $p(D_{I',J}) = D_{I,J}$.

{\it Case 2}: $I = (0,\dots,n)$ and $J = (0,\dots,n-1,n+1)$. Let $J' = (0,\dots,n)$, so $D_{I,J'} = D_0$, and let $$p = (-1)^{a}\frac{(n+1)!}{(a+n)!} \  \Omega_{a,0}\circ_{a-1}.$$ As above, $p(D_{I,J'}) = D_{I,J}$. Hence property $P_{n,k}$ holds for $k=(n+1)^2 + 1$.

\smallskip

{\bf Step 2}: By Step 1, $P_{1,5}$ holds, so assume inductively that $P_{1,k-1}$ holds. Let $D_{I,J}\in U_1[k]$. 

{\it Case 1}: $I = (i_0,i_1)\neq (0,1)$ and $i_1 - i_0 >1$. For any $b>j_1 + 1$, let $I' = (i_0,i_1-1)$ and let $$p = (-1)^{b+1} \frac{1}{(b+i_1-1)!} \ \Omega_{i_1, b}\circ_{i_1+b-1}.$$ By (\ref{actionp}), (\ref{paraction}), and (\ref{actioni}), we have $p(D_{I',J}) =D_{I,J}$.

{\it Case 2}: $I = (i_0,i_1)\neq (0,1)$ and $i_1 - i_0 = 1$. Then $i_0>0$, so we can take $I' = (i_0-1,i_1)$. Let $D_{I',J}$ be the corresponding element of $U_1[k-1]$. Also, let $I'' = (i_0-1,i_1+1)$ and let $D_{I'',J}$ be the corresponding element of $U_1[k]$. For any $b>j_1+1$, let $p_1 = \Omega_{0,b}\circ_{b-1}$ and $p_2 = \Omega_{1,b}\circ_b$. We have 
$$ p_1(D_{I',J}) = (-1)^{b+1} \frac{(b+i_0-1)!}{i_0!} \ D_{I,J} + (-1)^{b+1} \frac{(b+i_1)!}{(i_1+1)!}\ D_{I'',J},$$
$$ p_2(D_{I',J}) = (-1)^{b+1} \frac{(b+i_0-1)!}{(i_0-1)!} \ D_{I,J} + (-1)^{b+1} \frac{(b+i_1)!}{(i_1)!}\ D_{I'',J}.$$
Since the vectors $\big( \frac{(b+i_0-1)!}{i_0!}, \frac{(b+i_1)!}{(i_1+1)!}\big)$ and $\big(  \frac{(b+i_0-1)!}{(i_0-1)!}, \frac{(b+i_1)!}{(i_1)!}\big)$ in $\mathbb{C}^2$ are linearly independent, we can find a suitable linear combination $p = c_1 p_1 + c_2 p_2$ such that $p(D_{I',J}) =D_{I,J}$.

{\it Case 3}: $I = (0,1)$, $J = (j_0,j_1)\neq (0,1)$, and $j_1-j_0>1$. The argument is similar to Case 1 with the roles of $I$ and $J$ reversed.

{\it Case 4}: $I = (0,1)$, $J = (j_0,j_1)\neq (0,1)$, and $j_1-j_0 = 1$. This is similar to Case 2 with the roles of $I$ and $J$ reversed.

\smallskip

{\bf Step 3}: Assume that property $P_{m,k}$ holds for each pair $(m,k)$ with $m<n$ and $k\geq (m+1)^2$. We know from Step 1 that $P_{n,(n+1)^2+1}$ holds, so we may assume inductively that $P_{n,k-1}$ holds. Let $D_{I,J}\in U_n[k]$, and assume first that $I \neq (0,\dots,n)$. 

{\it Case 1}: $I \neq (0,\dots,n)$, and $i_0\neq 0$. Let $I'=(i_0-1,i_1,\dots,i_n)$, and let $D_{I',J}$ be the corresponding element of $U_n[k-1]$. For $b>j_n+1$, let $$p =  (-1)^{b+1} \frac{i_0!}{(b+i_0-1)!}\ \Omega_{i_0,b}\circ_{i_0+b-1}.$$ By (\ref{actionp}), (\ref{paraction}), and (\ref{actioni}), we have $p (D_{I',J}) = D_{I,J} + \sum_{r=1}^n \lambda_r D_{I^r,J}$, where $I^r = (i_0-1,i_1,\dots, i_{r-1}, i_r + 1,i_{r+1},\dots,i_n)$. For each $r=1,\dots,n$, let $$K^r = (i_1,\dots, i_{r-1}, i_r + 1, i_{r+1},\dots , i_n),$$ which is the list of length $n$ obtained from $I^r$ by removing the first entry $i_0-1$. Similarly, let $J' = (j_1,\dots,j_n)$, and let $D_{K^r,J'}$ be the corresponding element of $U_{n-1}$. By inductive assumption, for each $r=1,\dots,n$ there exist an integer $s_r> 0$ and a collection of lists $L^{r,1},\dots, L^{r,s_r}$ of length $n$, together with elements $p_{r,1},\dots, p_{r,s_r}\in \mathcal{P}$, such that 
$$ \sum_{t=1}^{s_r} p_{r,t}(D_{L^{r,t},J'}) = D_{K^r,J'}.$$
Moreover, we may assume that each of the lists $L^{r,t}$ has the property that the first term $l^{r,t}_1$ satisfies $0\leq  i^r_1-l^{r,t}_1 \leq 1$. Furthermore, we may assume that each $p_{r,t}$ is a linear combination of elements of the form $\Omega_{a,b}\circ_{a+b-1}$ with $a\geq i_1$ and $b>j_n+1$. 

Let $M^{r,t}$ be the list of length $n+1$ given by $(i_0-1,l^{r,t}_1,\dots,l^{r,t}_n)$, and let $D_{M^{r,t},J}$ be the corresponding element of $U_n$. Since $a-(i_0-1)\geq 2$, the operators $\Omega_{a,b}\circ_{a+b-1}$ do not affect the first index $i_0-1$ of $M^{r,t}$, so $\sum_{t=1}^{s_r} p_{r,t}(D_{M^{r,t},J}) = D_{I^r,J}$. Hence $$p(D_{I',J}) - \sum_{r=1}^n \sum_{t=1}^{s_r} \lambda_r p_{r,t}(D_{M^{r,t},J}) = D_{I,J}. $$

{\it Case 2}: $I \neq (0,\dots,n)$, and $i_0 = 0$. Let $r$ be the minimal integer for which $i_r>r$. By assumption, we have $1\leq r\leq n$. Hence $I = (0,\dots, r-1, i_r,\dots,i_n)$. Since $i_r - i_{r-1} = i_r - (r-1)\geq 2$, it follows that operators of the form $\Omega_{r+1,b}\circ_{r+b}$ only act on the last $n-r+1$ terms of $I$.

Let $I' = (i_r,\dots,i_n)$ and $J' = (j_r\dots,j_n)$ be the lists of length $n-r+1$ obtained from $I$ and $J$, respectively, by deleting the first $r$ terms, and let $D_{I',J'}\in U_{n-r}$ be the corresponding element. By inductive assumption there exist elements $D_{I^1,J'},\dots, D_{I^s,J'}\in U_{n-r}$ and elements $p_1,\dots,p_s\in\mathcal{P}$ such that \begin{equation}\label{stepiii} \sum_{t=1}^s p_t(D_{I^t,J'}) = D_{I',J'}.\end{equation} Moreover, we may assume that for $t=1,\dots,s$, $I^t = (i^t_r,\dots, i^t_n)$ satisfies $0\leq i_r - i^t_{r} \leq 1$, and that each $p_t$ is a linear combination of elements of the form $\Omega_{a,b}\circ_{a+b-1}$ where $a\geq i_r \geq r+1$ and $b>j_n+1$. 
Let $K^t = (0,\dots,r-1,i^t_r,\dots, i^t_n)$, and let $D_{K^t,J}$ be the corresponding element of $U_n$. It follows from (\ref{actionp}), (\ref{paraction}), (\ref{actioni}), and (\ref{stepiii}) that $$\sum_{t=1}^s p_t(D_{K^t,J}) = D_{I,J}.$$
 
{\it Case 3}:  $I = (0,\dots,n)$, $J\neq (0,\dots,n)$, and $j_0\neq 0$. The argument is the same as the proof of Case 1, with the roles of $I$ and $J$ reversed.

{\it Case 4}: $I = (0,\dots,n)$, $J\neq (0,\dots,n)$, and $j_0 = 0$. This is the same as Case 2, with the roles of $I$ and $J$ reversed. 

This completes the proof that property $P_{n,k}$ holds for all $n\geq 1$ and $k>(n+1)^2$. \end{proof}

\begin{remark}Specializing Theorem \ref{uniquesv} to the case $n=1$ proves the conjecture of Wang \cite{WIII} that all normally ordered polynomial relations in $\mathcal{W}_{1+\infty,-1}$ among the generators $j^0,j^1,j^2$ are consequences of a single relation.\end{remark}
There is a convenient way to think about the vertex operators $D_{I,J}\in \mathcal{M}_{-n}$ which is suggested by the proof of Lemma \ref{idealrecon}. Given a homogeneous polynomial $$p\in gr(\mathcal{M}_{-n})\cong \mathbb{C}[Q_{a,b}|\ a,b\geq 0]$$ of degree $k$ in the variables $Q_{a,b}$, a {\it normal ordering} of $p$ will be a choice of normally ordered polynomial $P\in (\mathcal{M}_{-n})_{(2k)}$, obtained by replacing $Q_{a,b}$ by $\Omega_{a,b}$, and replacing ordinary products with iterated Wick products of the form (\ref{iteratedwick}). Of course $P$ is not unique, but for any choice of $P$ we have $\phi_{2k}(P) = p$, where $\phi_{2k}: (\mathcal{M}_{-n})_{(2k)} \rightarrow (\mathcal{M}_{-n})_{(2k)} /(\mathcal{M}_{-n})_{(2k-1)} \subset gr(\mathcal{M}_{-n})$ is the usual projection. For the rest of this paper, $D^{2k}$, $E^{2k}$, $F^{2k}$, etc., will always denote elements of $(\mathcal{M}_{-n})_{(2k)}$ which are homogeneous, normally ordered polynomials of degree $k$ in the vertex operators $\Omega_{a,b}$.

Let $D_{I,J}^{2n+2}\in (\mathcal{M}_{-n})_{(2n+2)}$ be some normal ordering of $d_{I,J}$. Then $$\pi_{-n}(D_{I,J}^{2n+2}) \in (\mathcal{W}_{1+\infty,-n})_{(2n)},$$ where $\pi_{-n}:\mathcal{M}_{-n}\rightarrow\mathcal{W}_{1+\infty,-n}$ is the projection, and $\phi_{2n}(\pi_{-n}(D_{I,J}^{2n+2})) \in \text{gr}(\mathcal{W}_{1+\infty,-n}) \cong \mathbb{C}[q_{a,b}] / I_n$ can be expressed uniquely as a polynomial of degree $n$ in the variables $q_{a,b}$. Choose some normal ordering of the corresponding polynomial in the variables $\Omega_{a,b}$, and call this vertex operator $-D^{2n}_{I,J}$. Then $D^{2n+2}_{I,J} + D^{2n}_{I,J}$ has the property that $\pi_{-n}(D^{2n+2}_{I,J} + D^{2n}_{I,J})\in (\mathcal{W}_{1+\infty,-n})_{(2n-2)}.$ Continuing this process, we arrive at a vertex operator $D^{2n+2}_{I,J} + D^{2n}_{I,J}  + \cdots + D^2_{I,J}$ in the kernel of $\pi_{-n}$. We must have 
\begin{equation}\label{decompofd} D_{I,J} = \sum_{k=1}^{n+1}D^{2k}_{I,J},\end{equation} since $D_{I,J}$ is uniquely characterized by (\ref{uniquedij}).

In this decomposition, the term $D^2_{I,J}$ lies in the space $A_m$ spanned by $\{\Omega_{a,b}|\ a+b=m\}$, for $m = |I| + |J| + n$. Recall that $A_m  = \partial A_{m-1} \oplus \langle J^m \rangle$, and let $pr_m: A_m\rightarrow \langle J^m\rangle$ be the projection onto the second term. Define the {\it remainder} \begin{equation}\label{defofrij} R_{I,J} = pr_m(D^2_{I,J}).\end{equation} We will see that $R_{I,J}$ is independent of the choice of decomposition (\ref{decompofd}). 

\begin{lemma} \label{circzero} For all $j,k,l,m\geq 0$, $\Omega_{j,k}\circ_0\Omega_{l,m}$ is a total derivative. In particular, $\Omega_{j,k}\circ_0\Omega_{l,m}$ lies in the subspace $\partial(A_{s-1}) \subset A_s$, where $s = j+k +l+m$. Moreover, $\Omega_{j,k}\circ_1 \partial \Omega_{l,m}$ and $\Omega_{j,k}\circ_2 \partial^2 \Omega_{l,m}$ also lie in $\partial (A_{s-1})$.\end{lemma} 

\begin{proof} A calculation shows that for all $k,l\geq 0$, \begin{equation}\label{compi} J^k\circ_0 J^l = \sum_{i=0}^{k-1} (-1)^{i+1} \partial \Omega_{k+l-i-1,i}.\end{equation} We can write $\Omega_{j,k} = c_1 J^{j+k} + c_2 \partial \nu$, and $\Omega_{l,m} = d_1 J^{l+m} + d_2 \partial \mu$, for constants $c_1,c_2,d_1,d_2$ and vertex operators $\nu\in A_{j+k-1}$, $\mu\in A_{l+m-1}$. The first statement follows from (\ref{compi}) and the fact that $\partial \nu \circ_0 \Omega_{l,m} = 0$ and $J^{j+k} \circ_0 (\partial \mu) = \partial (J^{j+k} \circ_0 \mu) \subset \partial A_{s-1}$. The second statement then follows from the calculation $$\Omega_{j,k} \circ_1 \partial\Omega_{l,m} = \partial ( \Omega_{j,k} \circ_1 \Omega_{l,m}) - \partial \Omega_{j,k} \circ_1 \Omega_{l,m} = \partial ( \Omega_{j,k} \circ_1 \Omega_{l,m}) +  \Omega_{j,k} \circ_0 \Omega_{l,m},$$ and the third statement follows from $$\Omega_{j,k} \circ_2 \partial^2\Omega_{l,m} = \partial ( \Omega_{j,k} \circ_2 \partial \Omega_{l,m}) - \partial \Omega_{j,k} \circ_2 \partial \Omega_{l,m} = \partial ( \Omega_{j,k} \circ_2 \partial \Omega_{l,m}) +2 \Omega_{j,k} \circ_1 \partial \Omega_{l,m}. $$ \end{proof}

\begin{lemma} Let $E\in (\mathcal{M}_{-n})_{(2m)}$ be a vertex operator of degree $2m$, and choose a decomposition \begin{equation}\label{dcnew} E = \sum_{k=1}^m E^{2k},\end{equation} where $E^{2k}$ is a homogeneous, normally ordered polynomial of degree $k$ in the variables $\Omega_{a,b}$. If $E = \sum_{k=1}^m F^{2k}$ is any rearrangement of (\ref{dcnew}), i.e., another decomposition of the same form, then $E^{2}-F^2$ is a second derivative.\end{lemma}

\begin{proof} Let $\mu =\ :a_1 \cdots a_m:$ be a normally ordered {\it monomial} in $(\mathcal{M}_{-n})_{(2m)}$, where each $a_i$ is one of the generators $\Omega_{a,b}$. Let $\tilde{\mu} =\ :a_{i_1} \cdots a_{i_m}:$, where $(i_1,\dots, i_m)$ is some permutation of $(1,\dots,m)$. It suffices to show that any decomposition
$\mu - \tilde{\mu} = \sum_{k=1}^{m-1} E^{2k}$ of the difference $\mu - \tilde{\mu}\in (\mathcal{M}_{-n})_{(2m-2)}$, has the property that $E^2$ is a second derivative.

To prove this statement, we proceed by induction on $m$. For $m=1$, there is nothing to prove since $\mu - \tilde{\mu} =0$. For $m=2$, and $\mu =\ :\Omega_{a,b}\Omega_{c,d}:$, we have \begin{equation}\label{reari} \mu - \tilde{\mu} =\ :\Omega_{a,b}\Omega_{c,d}: - : \Omega_{c,d}\Omega_{a,b}: \ = \sum_{i\geq 0} \frac{(-1)^i}{(i+1)!}\partial^{i+1}(\Omega_{a,b}\circ_{i} \Omega_{c,d}),\end{equation} by (\ref{vaidii}). Since $\Omega_{a,b}\circ_{0} \Omega_{c,d}$ is already a total derivative by Lemma \ref{circzero}, it follows that $\mu -\tilde{\mu}$ is a second derivative, as claimed.

Next, we assume the result for $r\leq m-1$. Since the permutation group on $m$ letters is generated by the transpositions $(i,i+1)$ for $i=1,\dots,m-1$, we may assume without loss of generality that $$\tilde{\mu} = \ :a_1 \cdots a_{i-1} a_{i+1} a_i a_{i+2} \cdots a_m:.$$ If $i>1$, we have $\mu - \tilde{\mu} = \ :a_1 \cdots a_{i-1} f:$, where $f =\ :a_i \cdots a_m: - : a_{i+1}a_i a_{i+2} \cdots a_m:$, which lies in $(\mathcal{M}_{-n})_{(2m-2i)}$. Since each term of $f$ has degree at least $2$, it follows that $\mu - \tilde{\mu}$ can be expressed in the form $\sum_{k=i}^{m-1}E^{2k}$. Since $i>1$, there is no term of degree 2. Given any rearrangement $\mu - \tilde{\mu}=\sum_{k=1}^{m-1}F^{2k}$, it follows from our inductive hypothesis that the term $F^2$ is a second derivative.

Suppose next that $i=1$, so that $\tilde{\mu} = \ :a_2 a_1 a_3 \cdots a_m:$. Define $$\nu = \ :(:a_1 a_2:) a_3 \cdots a_m:\ ,\ \ \ \ \ \ \tilde{\nu} = \ :(:a_2 a_1:) a_3 \cdots a_m:\ ,$$ and note that $\nu  - \tilde{\nu} = \ : (:a_1 a_2: - : a_2 a_1: )f:$, where $f = \ :a_3 \cdots a_m:$. By (\ref{reari}), $:a_1 a_2:  - : a_2 a_1:$ is homogeneous of degree 2, so $\nu  - \tilde{\nu}$ is a linear combination of monomials of degree $2m-2$. By inductive assumption, any rearrangement $\nu - \tilde{\nu} =\sum_{k=1}^{m-1}F^{2k}$ has the property that $F^2$ is a second derivative.

Next, by (\ref{vaidi}), we have \begin{equation}\label{diffi} \mu - \nu = - \sum_{k\geq 0} \frac{1}{(k+1)!} \bigg( :(\partial^{k+1} a_1) (a_2 \circ_k f): + : (\partial^{k+1} a_2)(a_1\circ_k f):\bigg).\end{equation} Since the operators $\circ_k$ for $k\geq 0$ are homogeneous of degree $-2$, each term appearing in (\ref{diffi}) has degree at most $2m-2$. Moreover, $$deg \big(:(\partial^{k+1} a_1) (a_2\circ_k f):\big) = 2+ deg(a_2 \circ_k f),\ \ \ deg \big(:(\partial^{k+1} a_2) (a_1\circ_k f):\big) = 2+ deg(a_1 \circ_k f),$$ so the only way to obtain terms of degree $2$ is for $a_2\circ_k f$ or $a_1\circ_k f$ to be a scalar. This can only happen if $k>0$, in which case we obtain either $\partial^{k+1} a_1$ or $\partial^{k+1} a_2$, which are second derivatives. Finally, by inductive assumption, any rearrangement of $\mu - \nu$ can contain only second derivatives in degree 2.

Similarly, $\tilde{\mu} - \tilde{\nu}$ has degree at most $2m-2$, and any rearrangement of $\tilde{\mu} - \tilde{\nu}$ can only contain second derivatives in degree 2. Since $\mu - \tilde{\mu} = (\mu - \nu) + (\nu - \tilde{\nu}) + (\tilde{\nu} - \tilde{\mu})$, the claim follows. \end{proof}

\begin{cor} \label{uniquenessofr} Fix $n\geq 1$. Given $D_{I,J}\in U_n$, suppose that $D_{I,J} = \sum_{k=1}^{n+1} D^{2k}_{I,J}$ and $D_{I,J} = \sum_{k=1}^{n+1} \tilde{D}^{2k}_{I,J}$ are two different decompositions of $D_{I,J}$ of the form (\ref{decompofd}). Then $$D^2_{I,J} - \tilde{D}^2_{I,J} \in \partial^2 (A_{s-2}),$$ where $s = |I|+|J| + n$. In particular, $R_{I,J}$ is independent of the choice of decomposition of $D_{I,J}$.\end{cor}

\begin{lemma} Fix $n\geq 1$, and let $R_0$ denote the remainder of the element $D_0$ given by (\ref{defodo}). The condition $R_0 \neq 0$ is equivalent to the existence of a decoupling relation of the form $j^l = P(j^0,\dots, j^{l-1})$ in $\mathcal{W}_{1+\infty,-n}$, for $l = n^2 + 2n$.\end{lemma}

\begin{proof} Let $D_0 =  \sum_{k=1}^{n+1}D^{2k}_0$ be a decomposition of $D_0$ of the form (\ref{decompofd}). If $R_0\neq 0$, we have $D^2_0 = \lambda J^l + \partial \omega$ for some $\lambda \neq 0$ and some $\omega\in A_{l-1}$. Applying the projection $\pi_{-n}:\mathcal{M}_{-n}\rightarrow \mathcal{W}_{1+\infty,-n}$, since $\pi_{-n}(D_0)=0$ we obtain \begin{equation}\label{deca} j^l = -\frac{1}{\lambda}\big( \partial \pi_{-n}(\omega) + \sum_{k=2}^{n+1} \pi_{-n}(D^{2k}_0) \big),\end{equation} which is a decoupling relation of the desired form. The converse follows from the fact that $D_0$ is the unique element of $\mathcal{I}_{-n}$ of weight $(n+1)^2$. \end{proof}

In the cases $n=1$ and $n=2$, we can give explicit formulas for $D_0$ and $R_0$. These calculations were performed using Kris Thielemann's OPE package \cite{T}. For $n=1$, let $$D^4_0 = \ :\Omega_{0,0}\Omega_{1,1}:-:\Omega_{0,1}\Omega_{1,0}:,\ \ \ \ \ \ D^2_0 =  \frac{1}{2} \partial \Omega_{1,1} +  \frac{1}{3} J^3.$$ We have $\pi_{-1}(D^4_0+D^2_0) = 0$, so we have $$D_0 = D^4_0 + D^2_0,\ \ \ \ \ \ R_0 = \frac{1}{3} J^3.$$
Similarly, for $n=2$, let $$D^6_0 = \ : \Omega_{0,0} \Omega_{1,1} \Omega_{2,2}: - : \Omega_{0,0} \Omega_{1,2} \Omega_{2,1}: - : \Omega_{0,1}\Omega_{1,0} \Omega_{2,2}:$$ $$+ \ :\Omega_{0,1} \Omega_{1,2} \Omega_{2,0}: - : \Omega_{0,2}\Omega_{1,1}\Omega_{2,0}: + :\Omega_{0,2}\Omega_{1,0}\Omega_{2,1}: , $$  $$D^4_0 = - 
 1/12 :\Omega_{0,0}\Omega_{2,5}:\  + 1/4 :\Omega_{0,1}\Omega_{2,4}:\  - 1/6 :\Omega_{0,2}\Omega_{2,3}: \  - 1/3 :\Omega_{0,1}\Omega_{4,2}:  $$ $$+ 1/3 :\Omega_{0,2}\Omega_{4,1}:\  - 1/4: \Omega_{0,0}\Omega_{5,2}: \  + 
 1/4 :\Omega_{0,2}\Omega_{5,0}:\  - 1/5 :\Omega_{0,0}\Omega_{6,1}: \ + 1/5 :\Omega_{0,1}\Omega_{6,0}: $$ $$+ 
 2/3:\Omega_{1,1} \Omega_{2,3}: \ - 2/3 : \Omega_{1,3} \Omega_{2,1}:\  - 1/6 : \Omega_{1,0} \Omega_{2,4}:\  + 
 1/6 :\Omega_{1,4} \Omega_{2,0}: \ + 1/3 :\Omega_{1,1} \Omega_{3,2}: $$ $$- 1/3 :\Omega_{1,2} \Omega_{3,1}: \ +  1/4 :\Omega_{1,0} \Omega_{4,2}: \  - 1/4 :\Omega_{1,2} \Omega_{4,0}:\  + 1/5 :\Omega_{1,0} \Omega_{5,1}:\  - 
 1/5 :\Omega_{1,1} \Omega_{5,0}: $$ $$- 1/3 :\Omega_{2,0} \Omega_{3,2}:\  + 1/3 :\Omega_{2,2} \Omega_{3,0}: \ - 
 1/4 :\Omega_{2,0} \Omega_{4,1}: \ + 1/4  :\Omega_{2,1} \Omega_{4,0}: , $$ and 
 $$D^2_0 =  -\frac{1}{30} \partial \Omega_{2,5} - \frac{7}{180}  \partial \Omega_{3,4} - \frac{1}{10} \partial \Omega_{4,3} - \frac{1}{10}  \partial \Omega_{5,2}  + \frac{7}{180} \partial \Omega_{7,0}  - \frac{1}{120} J^8.$$ We have $D_0 = D^6_0+ D^4_{0} + D^2_{0}$, so that $R_{0} = -\frac{1}{120} J^8$. For $n>2$, it is difficult to find explicit formulas for $D_0$ and $R_0$, but we will show that $R_0\neq 0$.

\begin{lemma} \label{secder} Fix $n\geq 1$, and suppose that $D_{I,J}\in U_n$ has the property that $R_{I,J} = 0$. Then for any decomposition $D_{I,J} = \sum_{k=1}^{n+1}D^{2k}_{I,J}$ of the form (\ref{decompofd}), the term $D_{I,J}^2$ is a second derivative. In particular, $D_{I,J}^2 \in \partial^2(A_{s-2})$ where $s= |I|+|J|+n$.\end{lemma}

\begin{proof} Let $I = (i_0,\dots,i_n)$ and $J= (j_0,\dots,j_n)$. By (\ref{paraction}) and (\ref{actioni}), $D_{I,J}$ is an eigenvector of $J^2\circ_2$ with eigenvalue $\lambda = - \sum_{r=0}^n i_r(i_r-1) + \sum_{r=0}^{n} (j_r + 1)(j_r + 2)$, so \begin{equation}\label{jeigen} J^2\circ_2 \big(\sum_{k=1}^{n+1}D^{2k}_{I,J}) = \lambda \big(\sum_{k=1}^{n+1}D^{2k}_{I,J}\big).\end{equation}
On the other hand, we can compute $J^2\circ_2 \big(\sum_{k=1}^{n+1}D^{2k}_{I,J}\big)$ using (\ref{vaidiii}). For $k>2$, it follows from (\ref{vaidiii}) that $J^2\circ_2 D^{2k}_{I,J}$ can be expressed in the form $E^{2k} + E^{2k-2}+E^{2k-4}$. For $k>3$, any rearrangement of $J^2\circ_2 D^{2k}_{I,J}$ will only contain second derivatives in degree $2$, so we only need to consider the contribution from $J^2\circ_2 D^{2k}_{I,J}$ for $k=1$, $k=2$, and $k=3$.

For $k=3$, $D^6_{I,J}$ is a sum of terms of the form $:\Omega_{a,b} \Omega_{c,d}\Omega_{e,f}:$, and by (\ref{vaidiii}), the contribution of $J^2\circ_2 D^6_{I,J}$ in degree 2 will consist of terms of the form $$J^2\circ_0 (\Omega_{a,b} \circ_0 (\Omega_{c,d} \circ_0 \Omega_{e,f})),$$ which is a total derivative by Lemma \ref{circzero}.

For $k=2$, $D^4_{I,J}$ consists of a sum of terms of the form $:\Omega_{a,b}\Omega_{c,d}:$. By (\ref{vaidiii}), the contribution of $J^2\circ_2 D^4_{I,J}$ in degree 2 is a sum of terms of the form $(J^2\circ_0 \Omega_{a,b}) \circ_1 \Omega_{c,d}$ and $(J^2\circ_1 \Omega_{a,b}) \circ_0 \Omega_{c,d}$, which are all total derivatives, by Lemma \ref{circzero}. 

For $k=1$, $D^2_{I,J} = \sum_{i=1}^s c_i \partial^i J^{s-i}$, since $R_{I,J} = 0$. Suppose that $c_1\neq 0$. Using (\ref{actionp}), we calculate $$J^2\circ_2 \partial J^{s-1} = -(2s+2) J^s - (s^2-3s-4) \partial J^{s-1}.$$ By Lemma \ref{circzero}, $J^2\circ_2 ( \sum_{i=2}^s c_i \partial^i J^{s-i})$ is a total derivative, so $J^2 \circ_2 D^2_{I,J}\equiv -(2s+2)c_1 J^s $ modulo $\partial (A_{s-1})$. But $J^2\circ_2 D^2_{I,J} \equiv \lambda D^2_{I,J}$ modulo $\partial^2 (A_{s-2})$ by (\ref{jeigen}) and Corollary \ref{uniquenessofr}, which violates the fact that $R_{I,J} = 0$. \end{proof}

\begin{lemma} \label{nojo} Suppose that $D_{I,J}$ has the property that $R_{I,J} = 0$, and for some decomposition $D_{I,J} = \sum_{k=1}^{n+1} D^{2k}_{I,J}$, the term $D^4_{I,J}$ does not depend on $\Omega_{0,0}$. Then for any decomposition of $\Omega_{0,1}\circ_2 D_{I,J}$ of the form $\sum_{k=1}^{n+1}E^{2k}$, the term $E^2$ is a total derivative.\end{lemma}

\begin{proof} It suffices to find {\it some} decomposition of $\Omega_{0,1}\circ_2 D_{I,J}$ of the desired form; since the property of $E^2$ being a total derivative is stable under rearrangements, any other decomposition of $\Omega_{0,1}\circ_2 D_{I,J}$ will have this property as well. As in the proof of the preceding lemma, for $k\geq 3$, $\Omega_{0,1}\circ_2 D^{2k}_{I,J}$ can be expressed in the form $E^{2k} + E^{2k-2}+E^{2k-4}$. The same argument then shows that for $k\geq 3$, $\Omega_{0,1}\circ_2 D^{2k}_{I,J}$ can only contribute a second derivative in degree $2$.

For $k=2$, $D^4_{I,J}$ consists of a sum of terms of the form $:\Omega_{a,b}\Omega_{c,d}:$. By (\ref{vaidiii}), we have $\Omega_{0,1}\circ_2 (:\Omega_{a,b} \Omega_{c,d}:) =$ $$ :(\Omega_{0,1}\circ_2 \Omega_{a,b}) \Omega_{c,d}: + \ :\Omega_{a,b} (\Omega_{0,1}\circ_2 \Omega_{c,d}): \  + (\Omega_{0,1}\circ_1 \Omega_{a,b})\circ_0 \Omega_{c,d} + (\Omega_{0,1}\circ_0 \Omega_{a,b})\circ_1 \Omega_{c,d}.$$ Since $\Omega_{0,1}\circ_2$ lowers weight by 1, the only element of the form $\Omega_{a,b}$ for which $\Omega_{0,1}\circ_2 \Omega_{a,b}$ is a constant is $\Omega_{0,0}$. Since $D^4_{I,J}$ does not depend on $\Omega_{0,0}$, none of the terms $:(\Omega_{0,1}\circ_2 \Omega_{a,b}) \Omega_{c,d}:$ or $:\Omega_{a,b} (\Omega_{0,1}\circ_2 \Omega_{c,d}): $ appearing in $\Omega_{0,1}\circ_2 D^4_{I,J}$ will have degree 2. Moreover, each term of the form  $(\Omega_{0,1}\circ_1 \Omega_{a,b})\circ_0 \Omega_{c,d} + (\Omega_{0,1}\circ_0 \Omega_{a,b})\circ_1 \Omega_{c,d}$ appearing in $\Omega_{0,1} \circ_2 D^4_{I,J}$ is a total derivative, by Lemma \ref{circzero}. It follows that the component of $\Omega_{0,1}\circ_2 D^4_{I,J}$ in degree 2 is a total derivative. 

Finally, for $k=1$, it follows from Lemma \ref{secder} that $D^2_{I,J}$ is a second derivative, since $R_{I,J}=0$. By Lemma \ref{circzero}, $\Omega_{0,1}\circ_2 D^2_{I,J}$ is then a total derivative. We conclude that $\Omega_{0,1}\circ_2 D_{I,J}$ can be expressed in the form $\sum_{k=1}^{n+1}E^{2k}$ where $E^{2}$ is a total derivative, as claimed. \end{proof}

\begin{lemma} \label{inductivestep} Fix $n\geq 1$, and assume that $D_0$ has the property that the remainder $R_0\neq 0$. Let $I = (1,\dots, n+1) = J$, and let $D_{I,J}$ be the corresponding element of $U_n$. Then the remainder $R_{I,J}\neq 0$ as well.\end{lemma}

\begin{proof} Choose a decomposition $D_{I,J}= \sum_{k=1}^{n+1} D_{I,J}^{2k}$ of the form (\ref{decompofd}). We will assume that $R_{I,J}= 0$, and obtain a contradiction by showing that this implies that $R_0 = 0$. Since $D_{I,J}^{2n+2}$ does not depend on the variables $\{\Omega_{a,0}|\ a\geq 0\}$, we may assume without loss of generality that for $k=1,\dots,n$, $D^{2k}_{I,J}$ also does not depend on $\{\Omega_{a,0}|\ a\geq 0\}$. In particular, $D^4_{I,J}$ does not depend on $\Omega_{0,0}$ so the hypothesis of Lemma \ref{nojo} is satisfied.

By Lemma \ref{secder}, we may assume that $D^2_{I,J}$ is a second derivative. We act on $D_{I,J}$ by the operator $\Omega_{0,1}\circ_2$, which is homogenous of weight $-1$. For $k>0$, $\Omega_{0,1}\circ_2 \Omega_{j,k}$ cannot contain a term of the form $\Omega_{a,0}$, by (\ref{paraction}) and (\ref{actioni}). Similarly, if $P$ is any normally ordered polynomial in the variables $\{\Omega_{j,k}|\ k>0\}$, $\Omega_{0,1}\circ_2 P$ can be expressed as a normally ordered polynomial in $\{\Omega_{j,k}|\ k>0\}$ as well. Using (\ref{paraction}) and (\ref{actioni}), we calculate $$\Omega_{0,1}\circ_2 D_{I,J} = 2 D_{K^1,J},\ \ \ \ \ \ \ \ K^1 = (0,2,3,\dots,n+1).$$ Since $D_{I,J}^4$ is independent of $\Omega_{0,0}$, it follows from Lemma \ref{nojo} that $D_{K^1,J}$ also satisfies $R_{K^1,J} = 0$. Lemma \ref{secder} then shows that $D^2_{K^1,J}$ is a second derivative. Finally, since $D_{I,J}$ is a normally ordered polynomial in $\{\Omega_{j,k}|\ k>0\}$, we can assume that $D_{K^1,J}$ has a decomposition
$$D_{K^1,J} = \sum_{k=1}^{n+1}D^{2k}_{K^1,J} $$ such that each term $D^{2k}_{K^1,J}$ does not depend on $\{\Omega_{a,0}|\ a\geq 0\}$. Repeating this argument $n$ times, we find that $$\big(\Omega_{0,1}\circ_2 \big)^n (D_{I,J}) = \lambda D_{K,J},$$ where $K = (0,1,\dots,n)$ and $\lambda = \prod_{i=1}^{n+1} (i)(i+1)$. Moreover, $D_{K,J}$ has a decomposition $$D_{K,J} = \sum_{k=1}^{n+1}D^{2k}_{K,J}$$ such that $D^2_{K,J}$ is a second derivative (so in particular $R_{K,J} = 0$), and $D^4_{K,J}$ does not depend on $\{\Omega_{a,0}|\ a\geq 0\}$. In particular, $D^4_{K,J}$ is independent of $\Omega_{0,0}$.

Next, we act on $D_{K,J}$ by $\Omega_{0,0}\circ_1$. We get $\Omega_{0,0}\circ_1 D_{K,J} = D_{K,L^1}$, where $L^1 = 0,2,\dots,n+1$. Moreover, since $D^4_{K,J}$ is independent of $\Omega_{0,0}$, we can conclude that $D^2_{K,L^1}$ is a second derivative. However, we can no longer conclude that $D^4_{K,L^1}$ does not depend on $\Omega_{0,0}$. So instead of using the operator $\Omega_{0,0}\circ_1$ to lower the weight at this stage, we use the operator $$f = \Omega_{1,0}\circ_2 + \Omega_{0,0}\circ_1.$$ Clearly $f$ lowers the weight by 1, and by the same argument as the proof of Lemma \ref{circzero}, $f(\partial^2 A_{s}) \subset \partial A_{s}$. Moreover, $f(\Omega_{0,0}) = 0$, so given any normally ordered monomial $D^4$ of the form $:\Omega_{a,b}\Omega_{c,d}:$, the component of $f(D^4)$ in degree $2$ will be a total derivative, as in the proof of Lemma \ref{nojo}.

Note that $f (D_{K,L^1})= -4 D_{K,L^2}$ where $L^2 = (0,1,3,\dots,n+1)$. Since $D_{K,L^1}^2$ is a second derivative, $f(D_{K,L^1}^2)$ is a total derivative, and since any rearrangement of $f(D^{2k}_{K,L^1})$ can only contribute a second derivative in degree 2 for $k>1$, we conclude that $R_{K,L^2}$ is zero. Applying Lemma \ref{secder} again, $D^2_{K,L^2}$ is then a second derivative. For $i=2,\dots,n$ let $L^i = (0,1,\dots,i-1,i+1,\dots,n+1)$. It is easy to check using (\ref{paraction}) and (\ref{actioni}) that
$$f(D_{K,L^i}) = -\big(i + 1)^2D_{K,L^{i+1}}.$$ Moreover, at each stage, $D^2_{K,L^i}$ is a second derivative, and in particular, $R_{K,L^i} = 0$. At the $n$th stage, we see that $f(D_{K,L^n}) = -(n+1)^2 D_0$, so we have $R_0=0$, which is a contradiction. \end{proof}

Recall that $\mathcal{S}(V)$ is a graded algebra with $\mathbb{Z}_{\geq 0}$ grading (\ref{grading}), which specifies a linear isomorphism $\mathcal{S}(V)\cong Sym\bigoplus_{k\geq 0} (V_k \oplus V^*_k)$. By (\ref{gradingw}), $\mathcal{W}_{1+\infty,-n}$ is a graded subalgebra of $\mathcal{S}(V)$, so we obtain an isomorphism of graded vector spaces \begin{equation}\label{linisomor} i_{-n}: \mathcal{W}_{1+\infty,-n} \rightarrow (Sym\bigoplus_{k\geq 0} (V_k \oplus V^*_k))^{GL_n}.\end{equation}

Let $p\in \big(Sym\bigoplus_{k\geq 0} (V_k \oplus V^*_k )\big)^{GL_n}$ be a homogeneous polynomial of degree $2d$, and let $f = (i_{-n})^{-1}(p)\in (\mathcal{W}_{1+\infty,-n})^{(2d)}$ be the corresponding homogeneous vertex operator. Let $F\in (\mathcal{M}_{-n})_{(2d)}$ be a vertex operator satisfying $\pi_{-n}(F) = f$, where $\pi_{-n}: \mathcal{M}_{-n}\rightarrow \mathcal{W}_{1+\infty,-n}$ is the usual projection. We can write $F = \sum_{k=1}^{d} F^{2k}$, where $F^{2k}$ is a normally ordered polynomial of degree $k$ in the vertex operators $\Omega_{a,b}$.

Next, let $\tilde{V}$ be the vector space $\mathbb{C}^{n+1}$, and let $$\tilde{q}_{a,b}\in \big(Sym \bigoplus_{k\geq 0}( \tilde{V}_k \oplus \tilde{V}^*_k) \big)^{GL_{n+1}}$$ be the generator given by (\ref{weylgenerators}). Here $\tilde{V}_k$ and $\tilde{V}^*_k$ are isomorphic to $\tilde{V}$ and $\tilde{V}^*_k$, respectively, for $k\geq 0$. Let $\tilde{p}$ be the polynomial of degree $2d$ obtained from $p$ by replacing each $q_{a,b}$ with $\tilde{q}_{a,b}$, and let $\tilde{f} = (i_{-n-1})^{-1} (\tilde{p}) \in (\mathcal{W}_{1+\infty,-n-1})^{(2d)}$ be the corresponding homogeneous vertex operator. Finally, let $\tilde{F}^{2k}\in \mathcal{M}_{-n-1}$ be the vertex operator obtained from $F^{2k}$ by replacing each $\Omega_{a,b}$ with the corresponding vertex operator $\tilde{\Omega}_{a,b}\in \mathcal{M}_{-n-1}$, and let $\tilde{F} = \sum_{i=1}^d \tilde{F}^{2k}$. 

\begin{lemma} \label{indepofn} We can choose $F$ such that $\pi_{-n-1}(\tilde{F}) = \tilde{f}$.\end{lemma}

\begin{proof} We may assume without loss of generality that $p$ is a monomial in the variables $q_{a,b}$. If $d = 1$, $p=q_{a,b}$ for some $a,b\geq 0$, and $f = \omega_{a,b}$. We can take $F = \Omega_{a,b}$, so the claim is obvious. We assume inductively that for monomials $p =q_{a_1,b_1} \cdots q_{a_r,b_r}$ for $r<d$, there is a vertex operator $F = \sum_{k=1}^r F^{2k}$ such that $\pi_{-n}(F) = f$ where $f = (i_{-n})^{-1}(p)$, such that $\pi_{-n-1}(\tilde{F}) = \tilde{f}$, and $F^{2r} =\ : \Omega_{a_1,b_1}\cdots \Omega_{a_r,b_r}:$.

Now let $p = q_{a_1,b_1} \cdots q_{a_d,b_d}$. By inductive assumption, there exists a vertex operator $G = \sum_{k=1}^{d-1} G^{2k} \in \mathcal{M}_{-n}$ such that $$G^{2d-2} = \ :\Omega_{a_2,b_2}\cdots \Omega_{a_d,b_d}:\ ,\ \ \ \ \ \pi_{-n}(G) = g,\ \ \ \ \ \pi_{-n-1}(\tilde{G}) = \tilde{g},$$ where $g = (i_{-n})^{-1}(q_{a_2,b_2}\cdots q_{a_d,b_d})$. Define a vertex operator $H\in \mathcal{M}_{-n}$ by $$H = \sum_{k=2}^d H^{2k},\ \ \ \ \  H^{2k} = \ :\Omega_{a_1,b_1} G^{2k-2}:\ .$$
Since $\pi_{-n}$ is a vertex algebra homomorphism, we have $$\pi_{-n}(H) = \pi_{-n} \big(:\Omega_{a_1,b_1} G:\big) = \ : \omega_{a_1,b_1} g:,$$ and using (\ref{vaidi}), we see that \begin{equation}\label{nodouble} :\omega_{a_1,b_1} g:\  = f + f',\end{equation} where $f'$ is homogeneous of degree $2d-2$. In fact, a computation using (\ref{vaidi}) shows that under the isomorphism (\ref{linisomor}), $f'$ corresponds to the polynomial
$$\sum_{r=2}^d \bigg( \frac{(-1)^{b_1+1}}{a_r + b_1 + 1} q_{a_1 + b_1 + a_r,b_r} + \frac{(-1)^{a_1}}{a_1+b_r + 1} q_{a_r,b_r + a_1 + b_1} \bigg) q_{a_2,b_2}\cdots \widehat{q_{a_r,b_r}} \cdots q_{a_d,b_d}.$$ In this notation, the symbol $\widehat{q_{a_r,b_r}}$ means that the factor $q_{a_r,b_r}$ has been omitted, so the above polynomial is homogeneous of degree $2d-2$. Since this formula is independent of $n$, it follows that $\pi_{-n-1} (\tilde{H}) = \tilde{f} + \tilde{f}'$. By inductive assumption, there is a vertex operator $A = \sum_{k=1}^{d-1}A^{2k} \in \mathcal{M}_{-n}$ such that $\pi_{-n}(A) = f'$ and $\pi_{-n-1}(\tilde{A}) = \tilde{f}'$. Finally, we define $F = \sum_{k=1}^d F^{2k}$ by $F^{2d} = H^{2d}$, and $F^{2k} = H^{2k} - A^{2k}$ for $1\leq k< d$. It is immediate that $F$ has the desired properties. \end{proof}

\begin{cor} \label{corhomo} Fix $n\geq 1$, and let $D_{I,J}\in U_{n}$. There exists a decomposition $D_{I,J} = \sum_{k=1}^{n+1} D^{2k}_{I,J}$ of the form (\ref{decompofd}) such that the corresponding vertex operator $$\tilde{D}_{I,J} = \sum_{k=1}^{n+1} \tilde{D}^{2k}_{I,J} \in \mathcal{M}_{-n-1}$$ has the property that $\pi_{-n-1}(\tilde{D}_{I,J})$ lies in the homogeneous subspace $(\mathcal{W}_{1+\infty,-n-1})^{(2n+2)}$ of degree $2n+2$.\end{cor} 

\begin{proof} For each monomial $\mu$ of degree $2n+2$ appearing in the polynomial $d_{I,J}$, let $f_{\mu} = (i_{-n})^{-1}(\mu)\in (\mathcal{W}_{1+\infty,-n})^{(2n+2)}$. By the preceding lemma, we may choose $F_{\mu} = \sum_{k=1}^{n+1} F_{\mu}^{2k} \in \mathcal{M}_{-n}$ such that $\pi_{-n}(F_{\mu}) = f_\mu$, and $\pi_{-n-1}(\tilde{F}_{\mu}) = \tilde{f}_{\mu}$. For each $k=1,\dots,n+1$, define $D_{I,J}^{2k} = \sum_{\mu} F_{\mu}^{2k}$ where the sum is over all monomials $\mu$ appearing in $d_{I,J}$, and take $D_{I,J} = \sum_{k=1}^{n+1} D^{2k}_{I,J}$. Clearly this is a decomposition of $D_{I,J}$ of the form (\ref{decompofd}), and we have $$\pi_{-n-1}(\tilde{D}_{I,J}) = \sum_{\mu} \pi_{-n-1} (\tilde{F}_{\mu}) = \sum_{\mu} \tilde{f}_{\mu},$$ which is homogeneous of degree $2n+2$. \end{proof}

Now we have assembled all the technical tools necessary to prove our main result.

\begin{thm} \label{decoupexist} For all $n\geq 1$, we have $R_0 \neq 0$. Hence there is a decoupling relation in $\mathcal{W}_{1+\infty,-n}$ of the form $j^l = P(j^0,\dots, j^{l-1})$, for $l = n^2 + 2n$.\end{thm} 

\begin{proof} This is well known for $n=1$ and we have already shown it for $n=2$, so we will assume inductively that it holds for $n-1$. The idea of the proof is to use our inductive assumption together with Lemma \ref{inductivestep} to construct a decomposition
$$D_0 = \sum_{k=1}^{n+1} D^{2k}_0$$
with the property that $D^4_0$ contains the term $: J^0 J^{l-1}:$ with nonzero coefficient, for $l=n^2+2n$. By (\ref{reari}), $:J^0 J^{l-1}: - :J^{l-1} J^0:$ is a second derivative, so we may assume that $:J^{l-1} J^0:$ does not appear in $D^4_0$. Suppose that such a decomposition exists, and that $R_0= 0$. First, (\ref{vaidiii}) shows that for $k>1$, $J^0\circ_1D_0^{2k}$ can be expressed in the form $E^{2k} + E^{2k-2}$, so for $k>2$, any rearrangement of $J^0\circ_1D_0^{2k}$ can only contribute a second derivative in degree $2$. Moreover, since $R_0 = 0$, $D_0^2$ is a total derivative, so by Lemma \ref{circzero}, $J^0\circ_1D_0^2$ is also a total derivative. Since $D^4_0$ contains the term $:J^0 J^{l-1}:$ with nonzero coefficient, $J^0\circ_1D^4_0$ will contain the term $J^{l-1}$ with nonzero coefficient, and this term cannot be canceled by any term coming from $J^0\circ_1 D^{2k}_0$ for $k\neq 2$. This contradicts the fact that $D_0$ is a singular vector.

Let $d_0\in gr(\mathcal{M}_{-n})\cong \mathbb{C}[Q_{a,b}|\ a,b\geq 0]$ denote the image of $D_0$ under the projection $\phi_{2n+2}: (\mathcal{M}_{-n})_{(2n+2)}\rightarrow (\mathcal{M}_{-n})_{(2n+2)}/(\mathcal{M}_{-n})_{(2n+1)}\subset gr(\mathcal{M}_{-n})$. By Theorem \ref{weylfft}, $$d_0 = \det  \left[\begin{matrix}Q_{0,0} & \cdots & Q_{0,n} \cr  \vdots  & & \vdots  \cr  Q_{n,0}  & \cdots & Q_{n,n} \end{matrix} \right],$$ so $d_{0}$ can be written in the form \begin{equation}\label{decdet} d_0 = Q_{0,0} d_{I,J} + d',\end{equation} where $I = (1,\dots,n) = J$, $d_{I,J}$ is the corresponding polynomial of degree $2n$, and $d'$ is a polynomial of degree $2n+2$ which does not depend on $Q_{0,0}$. Consider the vertex operator $D_{I,J}\in \mathcal{M}_{-n+1}$ corresponding to $d_{I,J}$, regarded now as an element of $gr(\mathcal{M}_{-n+1})$. By Corollary \ref{corhomo}, we may choose a decomposition
$D_{I,J} = \sum_{k=1}^n D_{I,J}^{2k}$ such that the corresponding vertex operator $\tilde{D}_{I,J} = \sum_{k=1}^n \tilde{D}_{I,J}^{2k} \in \mathcal{M}_{-n}$ has the property that $\pi_{-n}(\tilde{D}_{I,J})\in (\mathcal{W}_{1+\infty,-n})^{(2n)}$. Moreover, since $d_{I,J}$ does not depend on $Q_{0,0}$, we may assume that each term $D^{2k}_{I,J}$ appearing in $D_{I,J}$ is independent of $J^0 = \Omega_{0,0}$. We will use this decomposition of $\tilde{D}_{I,J}$ to create a decomposition of $D_{0}$ with the desired property. 

By our inductive assumption together with Lemma \ref{inductivestep}, $D^2_{I,J}$ contains the term $J^{l-1}$ with nonzero coefficient, and since $\tilde{D}^2_{I,J}$ is obtained from $D^2_{I,J}$ by replacing each $\Omega_{j,k}\in \mathcal{M}_{-n+1}$ with the corresponding element of $\mathcal{M}_{-n}$, it follows that $\tilde{D}^2_{I,J}$ contains $J^{l-1}$ (regarded now as an element of $\mathcal{M}_{-n}$), with nonzero coefficient as well. Consider the vertex operator $$: J^0  \tilde{D}_{I,J}:\  = \sum_{k=1}^n : J^0  \tilde{D}_{I,J}^{2k}:\ \in \mathcal{M}_{-n}.$$ Since $\tilde{D}^2_{I,J}$ has a nonzero term of the form $J^{l-1}$, $: J^0 \tilde{D}_{I,J}:$ has a nonzero term of the form $: J^0 J^{l-1}:$ in degree $4$. Applying the projection $\pi_{-n}:\mathcal{M}_{-n}\rightarrow \mathcal{W}_{1+\infty,-n}$, we have $$\pi_{-n}(:J^0  \tilde{D}_{I,J}:) = \ : j^0 f: ,$$ where $f = \pi_{-n} (\tilde{D}_{I,J})$. Since $f$ is homogeneous of degree $2n$, it follows from (\ref{vaidi}) that $:j^0 f: \  = g + g'$ where $g\in (\mathcal{W}_{1+\infty,-n})^{(2n+2)}$ and $g'\in (\mathcal{W}_{1+\infty,-n})^{(2n)}$. It is easy to see from (\ref{vaidi}) that under the isomorphism $i_{-n}$ given by (\ref{linisomor}), $i_{-n}(g')\in (Sym\bigoplus_{k\geq 0} (V_k \oplus V^*_k))^{GL_n}$ does not depend on $q_{0,0}$. Hence we can choose a vertex operator $E = \sum_{k=1}^n E^{2k}\in\mathcal{M}_{-n}$ such that $\pi_{-n}(E) = -g'$, and such that each $E^{2k}$ is independent of $J^0$. It follows that $\pi_{-n} ( :J^0  \tilde{D}_{I,J}: + E)\in (\mathcal{W}_{1+\infty,-n})^{(2n+2)}$.

Next, let $F^{2n+2}$ be a normal ordering of the polynomial $d'$ given by (\ref{decdet}). Since $d'$ is independent of $Q_{0,0}$, we may assume that $F^{2n+2}$ does not depend on $J^0$. Then $\pi_{-n}(F^{2n+2})$ will contain terms of lower degree in $\mathcal{W}_{1+\infty,-n}$, and we can find elements $F^{2k}\in \mathcal{M}_{-n}$ for $k=1,\dots,n$, such that $$\pi_{-n}(F)\in (\mathcal{W}_{1+\infty,-n})^{(2n+2)},\ \ \ \ \ \ F =  \sum_{k=1}^{n+1} F^{2k}.$$ Moreover, we may assume that each term $F^{2k}$ is independent of $J^0$. Finally, we define the decomposition \begin{equation}\label{finaldecomp} D_0 = \sum_{k=1}^{n+1} D^{2k}_0,\end{equation} where $D^2_0 = E^2 + F^2$, and $D^{2k}_0 = \ :J^0 \tilde{D}^{2k-2}_{I,J}:\ + E^{2k} + F^{2k}$ for $1<k\leq n+1$. Since $\pi_{-n}( :J^0 \tilde{D}_{I,J}:\  + E)$ and $\pi_{-n}(F)$ are both homogeneous of degree $2n+2$, and $$\phi_{2n+2}(:J^0 \tilde{D}_{I,J}:\ +E + F) = d_0\in gr(\mathcal{M}_{-n}),$$ it follows that (\ref{finaldecomp}) is indeed a decomposition of $D_0$, as claimed. Since $\tilde{D}^2_{I,J}$ contains the term $J^{l-1}$ with nonzero coefficient, and both $E^4$ and $F^4$ are independent of $J^0$, it follows that $D_0^4$ contains the term $: J^0 J^{l-1}:$ with nonzero coefficient, as desired. \end{proof}

\begin{thm} \label{stronggen} For $n\geq 1$, $\mathcal{W}_{1+\infty,-n}$ is strongly generated as a vertex algebra by $$\{j^0, j^1,\dots, j^{n^2+2n-1} \}.$$\end{thm}
\begin{proof} The decoupling relation $j^l= P(j^0,\dots,j^{l-1})$ for $l=n^2+2n$ given by Theorem \ref{decoupexist} is equivalent to the existence of an element $J^l - P(J^0,\dots,J^{l-1})\in \mathcal{I}_{-n}$. It suffices to show that for all $r>l$, there exists an element $J^r - Q_r(J^0,\dots,J^{l-1})\in \mathcal{I}_{-n}$, so we assume inductively that such an element exists for $r-1$.

Choose a decomposition $$Q_{r-1} = \sum_{k=1}^{d} Q^{2k}_{r-1},$$ where $Q_{r-1}^{2k}$ is a homogeneous normally ordered polynomial of degree $k$ in the vertex operators $J^0,\dots,J^{l-1}$ and their derivatives. In particular, $$Q^2_{r-1} = \sum_{i=0}^{l-1} c_i \partial^{r-i-1} J^{i},$$ for constants $c_0,\dots, c_{l-1}$.
We apply the operator $\Omega_{0,2}\circ_1\in \mathcal{P}$, which raises the weight by one. By (\ref{actionp}), we have $\Omega_{0,2}\circ_1 J^{r-1} = -(r+1)J^r$. Moreover, $\Omega_{0,2}\circ_1\big(\sum_{k=1}^{d} Q^{2k}_{r-1}\big)$ can be expressed in the form $\sum_{k=1}^{d} E^{2k}$ where each $E^{2k}$ is a normally ordered polynomial in $J^0,\dots,J^{l}$ and their derivatives. If $J^l$ or its derivatives appear in $E^{2k}$, we can use the element $J^l -P(J^0,\dots, J^{l-1})$ in $\mathcal{I}_{-n}$ to eliminate the variable $J^l$ and any of its derivatives, modulo $\mathcal{I}_{-n}$. Hence $\Omega_{0,2}\circ_1\big(\sum_{k=1}^{d} Q^{2k}_{r-1}\big)$ can be expressed modulo $\mathcal{I}_{-n}$ in the form $\sum_{k=1}^{d'} F^{2k}$, where $d'\geq d$, and $F^{2k}$ is a normally ordered polynomial in $J^0,\dots, J^{l-1}$ and their derivatives. It follows that $$-\frac{1}{r+1} \Omega_{0,2}\circ_1 \big(J^{r-1} - Q_{r-1}(J^0,\dots, J^{l-1})\big)$$ can be expressed as an element of $\mathcal{I}_{-n}$ of the desired form. \end{proof}

We remark that a similar strategy can be used to reprove the result from \cite{FKRW} that for $n\geq 1$, $\mathcal{M}_{n}$ has a unique singular vector $D$ of weight $n+1$, and $\mathcal{W}_{1+\infty,n}$ has a decoupling relation $j^n = P(j_0,\dots, j^{n-1})$. Recall that $\mathcal{W}_{1+\infty,n}$ can be realized as the invariant space $\mathcal{E}(V)^{GL_n}$, where $\mathcal{E}(V)$ is the $bc$-system associated to the vector space $V=\mathbb{C}^n$ \cite{FKRW}. The associated graded algebra $gr(\mathcal{E}(V))$ is $\bigwedge \bigoplus_{k\geq 0}(V_k\oplus V^*_k)$, and we have a linear isomorphism $$\mathcal{E}(V)^{GL_n} \cong \big(\bigwedge \bigoplus_{k\geq 0}(V_k\oplus V^*_k)\big)^{GL_n}.$$ There is a singular vector $D$ in $\mathcal{M}_{n}$ of weight $n+1$, which corresponds to the relation $(p_{0,0})^{n+1}\in \big(\bigwedge \bigoplus_{k\geq 0}(V_k\oplus V^*_k)\big)^{GL_n}$. Here $p_{0,0} = \sum_{i=1}^n x_{i,0}\wedge x'_{i,0}$, which is analogous to the corresponding element $q_{0,0}\in \big(Sym \bigoplus_{k\geq 0}(V_k\oplus V^*_k)\big)^{GL_n}$. We have a decomposition $D = \sum_{k=1}^{n+1} D^{2k}$, where $D^2$ is a linear combination of the vertex operators $\partial^i J^{n-i}$. An argument similar to the proof of Theorem \ref{decoupexist} shows that for all $n\geq 1$, the coefficient of $J^n$ in $D^2$ is nonzero, which yields a decoupling relation $j^n = P(j^0,\dots, j^{n-1})$ in $\mathcal{W}_{1+\infty,n}$. Finally, an argument analogous to the proof of Theorem \ref{uniquesv} shows that $D$ is the unique singular vector in $\mathcal{I}_{n}$.

A remaining question is whether the isomorphism $\mathcal{W}_{1+\infty,-1}\cong \mathcal{H}\otimes \mathcal{W}_{3,-2}$ given by Wang in \cite{WI} has an analogue for $n>1$. In other words, can $\mathcal{W}_{1+\infty,-n}$ be related to $\mathcal{W}(\mathfrak{g})$ for some Lie algebra $\mathfrak{g}$? We cannot answer this question at present, but there are a few things we can say. First, note that $j^0$ generates a copy of the Heisenberg algebra $\mathcal{H}$, and $j^0\circ_0$ acts by zero on $\mathcal{W}_{1+\infty,-n}$. Since $\mathcal{W}_{1+\infty,-n}$ decomposes as a direct sum of irreducible $\mathcal{H}$-modules, it follows that $\mathcal{W}_{1+\infty,-n} \cong \mathcal{H}\otimes \mathcal{A}$, where $\mathcal{A}$ is the commutant $Com(\mathcal{H},\mathcal{W}_{1+\infty,-n})$. Define vertex operators \begin{equation}\label{defofl} L = \frac{1}{2n}\bigg(:j^0 j^0:\ +n \partial j^0 - 2n j^1\bigg),\end{equation}
\begin{equation}\label{defofw} W = \ : j^0 j^0 j^0:\ +\frac{3n}{2} :j^0 \partial j^0:\  -3n :j^0  j^1:\  + \frac{n^2}{4} \partial^2 j^0 -\frac{3n^2}{2}  \partial j^1+ \frac{3n^2}{2} j^2.\end{equation}
Since $\mathcal{W}_{1+\infty,-n}$ is generated by $j^0$, $j^1$, and $j^2$, and we can use (\ref{defofl}) and (\ref{defofw}) to express $j^1$ and $j^2$ in terms of the vertex operators $j^0$, $L$ and $W$, it follows that $\{j^0,L,W\}$ is another generating set for $\mathcal{W}_{1+\infty,-n}$. A straightforward OPE calculation shows that both $L$ and $W$ commute with $j^0$, so that $\mathcal{A}$ is generated by $L$ and $W$ as a vertex algebra. Moreover, $L$ generates a Virasoro algebra with central charge $-n-1$ and $W$ is primary of conformal weight $3$. Thus we have proved
\begin{thm}$\mathcal{A}$ is a conformal vertex algebra with central charge $-n-1$.\end{thm} In the case $n=1$, an OPE calculation shows that $L$ and $W$ generate a copy of the Zamolodchikov $\mathcal{W}_3$ algebra with central charge $-2$, so we recover Wang's result. 

\section{The representation theory of $\mathcal{W}_{1+\infty,-n}$}
In \cite{WII}, Wang showed that the irreducible, highest-weight modules over $\mathcal{W}_{1+\infty,-1}$ correspond to the points on a certain complex algebraic variety of dimension 2. The key step was to compute the Zhu algebra of $\mathcal{W}_{1+\infty,-1}$. Given a vertex algebra $\mathcal{V}$ with weight grading $\mathcal{V} = \bigoplus_{n\in\mathbb{Z}} \mathcal{V}_n$, the {\it Zhu functor} attaches to $\mathcal{V}$ an associative algebra $A(\mathcal{V})$, together with a surjective linear map $\pi_{Zh}:\mathcal{V}\rightarrow A(\mathcal{V})$ \cite{Z}. For $a\in \mathcal{V}_{m}$ and $b\in\mathcal{V}$, define
\begin{equation}\label{defzhu} a*b = Res_z \bigg (a(z) \frac{(z+1)^{m}}{z}b\bigg),\end{equation} and extend $*$ by linearity to a bilinear operation $\mathcal{V}\otimes \mathcal{V}\rightarrow \mathcal{V}$. Let $O(\mathcal{V})$ denote the subspace of $\mathcal{V}$ spanned by elements of the form \begin{equation}\label{zhuideal} a\circ b = Res_z \bigg (a(z) \frac{(z+1)^{m}}{z^2}b\bigg)\end{equation} where $a\in \mathcal{V}_m$, and let $A(\mathcal{V})$ be the quotient $\mathcal{V}/O(\mathcal{V})$, with projection $\pi_{Zh}:\mathcal{V}\rightarrow A(\mathcal{V})$. Then $O(\mathcal{V})$ is a two-sided ideal in $\mathcal{V}$ under the product $*$, and $(A(\mathcal{V}),*)$ is a unital, associative algebra. The assignment $\mathcal{V}\mapsto A(\mathcal{V})$ is functorial, and if $\mathcal{I}$ is a vertex algebra ideal of $\mathcal{V}$, we have $A(\mathcal{V}/\mathcal{I})\cong A(\mathcal{V})/ I$, where $I = \pi_{Zh}(\mathcal{I})$.

The main application of the Zhu functor is to study the representation theory of $\mathcal{V}$. A $\mathbb{Z}_{\geq 0}$-graded module $M = \bigoplus_{n\geq 0} M_n$ over $\mathcal{V}$ is called {\it admissible} if for every $a\in\mathcal{V}_m$, $a(n) M_k \subset M_{m+k -n-1}$, for all $n\in\mathbb{Z}$. Given $a\in\mathcal{V}_m$, the Fourier mode $a(m-1)$ acts on each $M_k$. The subspace $M_0$ is then a module over $A(\mathcal{V})$ with action $[a]\mapsto a(m-1) \in End(M_0)$. In fact, $M\mapsto M_0$ provides a one-to-one correspondence between irreducible, admissible $\mathcal{V}$-modules and irreducible $A(\mathcal{V})$-modules. 

Let $\mathcal{V}$ be a vertex algebra which is strongly generated by a set of weight-homogeneous elements $\alpha_i$ of weights $w_i$, for $i$ in some index set $I$. Then $A(\mathcal{V})$ is generated by $\{ a_i = \pi_{Zh}(\alpha_i(z))|\ i\in I\}$. Moreover, $A(\mathcal{V})$ inherits a filtration (but not a grading) by weight, and the associated graded object $gr(A(\mathcal{V}))$ is a commutative algebra with generators $\{\bar{a}_i|\ i\in I\}$.
Given an element $f\in A(\mathcal{V})$ of weight at most $w$, let $\bar{f}\in gr(A(\mathcal{V}))$ denote the symbol of $f$, i.e., the image of $f$ in $gr(A(\mathcal{V}))[w]$.

Let $C_2(\mathcal{V})$ denote the the vector space spanned by elements $\{:(\partial \alpha) \beta:|\  \alpha,\beta\in \mathcal{V}\}$. It is well-known that $\mathcal{V} / C_2(\mathcal{V})$ is a commutative algebra. Moreover, the map $\mathcal{V} / C_2(\mathcal{V}) \rightarrow gr(A(\mathcal{V}))$ sending the coset of $\alpha\in \mathcal{V}$ to the symbol of $\pi_{Zh}(\alpha)$, is a surjective algebra homomorphism.

Now we consider the case of $\mathcal{M}_c$ and $\mathcal{W}_{1+\infty,c}$. For any $c\in \mathbb{C}$, $A(\mathcal{M}_c)$ is the polynomial algebra $\mathbb{C}[a^0, a^1,a^2,\cdots]$, where $a^l = \pi_{Zh}(J^l(z))$ \cite{FKRW}. Moreover, $A(\mathcal{W}_{1+\infty,c}) \cong  \mathbb{C}[a^0,a^1,a^2,\dots, ] / I_c$, where $I_c = \pi_{Zh}(\mathcal{I}_c)$, and we have a commutative diagram
\begin{equation}\label{commdiag} \begin{array}[c]{ccc}
\mathcal{M}_c &\stackrel{\pi_c}{\rightarrow}&\mathcal{W}_{1+\infty,c} \\
\downarrow\scriptstyle{\pi_{Zh}}&&\downarrow\scriptstyle{\pi_{Zh}}\\
A(\mathcal{M}_c) &\stackrel{A(\pi_c)}{\rightarrow}&A(\mathcal{W}_{1+\infty,c})
\end{array} .\end{equation} Since $A(\mathcal{W}_{1+\infty,c})$ is a commutative algebra, its irreducible modules are all one-dimensional.

If $c$ is an integer $n\geq 1$, it is known that $A(\mathcal{W}_{1+\infty,n})\cong \mathbb{C}[a^0,a^1,\dots, a^{n-1}]$ \cite{FKRW}. The irreducible $A(\mathcal{W}_{1+\infty,n})$-modules (and hence the irreducible, admissible $\mathcal{W}_{1+\infty,n}$-modules) then correspond to the points in $\mathbb{C}^n$. The situation is much more interesting in the case of negative integral central charge. For $n=1$, it was shown in \cite{WII} that $A(\mathcal{W}_{1+\infty,-1}) \cong \mathbb{C}[h,t,w] / I$, where $I$ is the ideal generated by $f(t,w) = w^2 -\frac{1}{9} t^2(8t+1)$. It follows that the irreducible, admissible $\mathcal{W}_{1+\infty,-1}$-modules are parametrized by the points on the variety $V(I)\subset \mathbb{C}^3$, which is just the product of an affine line with a rational curve. For $n>1$, it is immediate from Theorem \ref{stronggen} that $A(\mathcal{W}_{1+\infty,-n})$ is generated by $\{a^0, \dots, a^{n^2+2n-1}\}$. Hence $$A(\mathcal{W}_{1+\infty,-n})\cong \mathbb{C}[a^0,\dots, a^{n^2+2n-1}]/ I_{-n},$$ where $I_{-n}$ is now regarded as an ideal inside $\mathbb{C}[a^0,\dots, a^{n^2+2n-1}]$. Let $V(I_{-n})\subset \mathbb{C}^{n^2+2n}$ be the corresponding variety, which then parametrizes the irreducible, admissible modules over $\mathcal{W}_{1+\infty,-n}$.

\begin{thm} \label{subvariety} For all $n\geq 1$, $V(I_{-n})$ is a proper, closed subvariety of $\mathbb{C}^{n^2+2n}$.\end{thm}

\begin{proof} We need to construct a nontrivial relation among the generators $a^0,\dots,a^{n^2+2n-1}$ of $A(\mathcal{W}_{1+\infty,-n})$. Recall the vector space $U_n=  (\mathcal{M}_{-n})_{(2n+2)}\cap \mathcal{I}_{-n}$ given by (\ref{vsu}), whose component $U_n[k]$ of weight $k$ has a basis $\{D_{I,J}|\  |I|+ |J| + n+ 1 = k\}$. A basis for $U_n[(n+1)^2 + 2]$ consists of the following five elements:
$$D_{(0,\dots,n),(0,\dots,n-1,n+2)}, \ \ \ \ \ \ D_{(0,\dots,n),(0,\dots,n-2,n,n+1)},\ \ \ \ \ \ D_{(0,\dots,n-1,n+1),(0,\dots,n-1,n+1)},$$  $$D_{(0,\dots,n-2,n,n+1),(0,\dots,n)},\ \ \ \ \ \ D_{(0,\dots,n-1,n+2),(0,\dots,n)}.$$
 As in the proof of Lemma \ref{inductivestep}, let $f$ be the operator $\Omega_{1,0}\circ_2 + \Omega_{0,0}\circ_1$, which lowers the weight by one. We calculate $$f (D_{(0,\dots,n), (0,\dots,n-1,n+2)}) = -(n+2)^2 D_{(0,\dots,n), (0,\dots,n-1,n+1)},$$ $$f (D_{(0,\dots,n), (0,\dots,n-1,n+1)}) = -(n+1)^2 D_{(0,\dots,n),(0,\dots,n)},$$  $$f(D_{(0,\dots,n), (0,\dots,n-2,n,n+1)}) = -n^2 D_{(0,\dots,n), (0,\dots,n-1,n+1)}.$$ Since the remainder $R_0 = R_{(0,\dots,n),(0,\dots,n)}$ is nonzero by Theorem \ref{decoupexist}, it follows that $R_{(0,\dots,n), (0,\dots,n-1,n+2)}$ and $R_{(0,\dots,n), (0,\dots,n-2,n,n+1)}$ are both nonzero as well. Hence there exists a unique nontrivial linear combination \begin{equation}\label{defofe} E =  D_{(0,\dots,n), (0,\dots,n-1,n+2)} + \lambda D_{(0,\dots,n), (0,\dots,n-2,n,n+1)},\ \ \ \ \lambda \in\mathbb{C}\setminus \{0\},\end{equation} such that for any decomposition $E = \sum_{k=1}^{n+1} E^{2k}$, the term $E^2$ does not depend on $J^{n^2+2n+2}$, and hence is a total derivative.

Recall from the proof of Theorem \ref{decoupexist} that $D_0 = D_{(0,\dots,n),(0,\dots,n)}$ admits a decomposition $D_0 = \sum_{k=1}^{n+1} D_0^{2k}$ for which $D_0^4$ contains the term $:J^0 J^{n^2+2n-1}:$ with nonzero coefficient. By the same argument, there is a decomposition $$D_{(0,\dots,n),(0,\dots,n-1,n+1)} = \sum_{k=1}^{n+1} D^{2k}_{(0,\dots,n),(0,\dots,n-1,n+1)}$$ for which $D^4_{(0,\dots,n),(0,\dots,n-1,n+1)}$ contains $:J^0 J^{n^2+2n}:$ with nonzero coefficient. Similarly, there exist decompositions $$D_{(0,\dots,n),(0,\dots,n-1,n+2)} = \sum_{k=1}^{n+1} D^{2k}_{(0,\dots,n),(0,\dots,n-1,n+2)},$$ $$D_{(0,\dots,n),(0,\dots,n-2,n,n+1)} = \sum_{k=1}^{n+1} D^{2k}_{(0,\dots,n),(0,\dots,n-2,n,n+1)},$$ such that $D^4_{(0,\dots,n),(0,\dots,n-1,n+2)}$ and $D^4_{(0,\dots,n),(0,\dots,n-2,n,n+1)}$ both contain $:J^0 J^{n^2+2n+1}:$ with nonzero coefficient.

Let $E = \sum_{k=1}^{n+1} E^{2k}$ be a decomposition of $E$, where $E^{2k}$ is a homogeneous, normally ordered polynomial of degree $k$ in the variables $\Omega_{a,b}$. Under the linear change of variables (\ref{lincomb}), we may regard $E^{2k}$ as a polynomial in the variables $\partial^i J^l$, $i,l\geq 0$. Since $E$ has weight $(n+1)^2 +2$, and $E^2$ is a total derivative, $E$ only depends on $J^0,\dots,J^{n^2+2n+1}$ and their derivatives.

By weight considerations, for $k>3$, $E^{2k}$ only depends on $J^0,\dots,J^{n^2+2n-1}$ and their derivatives. We may also assume that $E^6$ only depends on $J^0,\dots,J^{n^2+2n-1}$ and their derivatives, since our decomposition can be chosen so that $:J^0J^0 J^{n^2+2n}:$ does not appear in $E^6$. The possible terms in $E^4$ which can contain either $J^{n^2+2n}$ or $J^{n^2+2n+1}$ are $:J^0 J^{n^2+2n+1}:$, $:J^1 J^{n^2+2n}:$, and $:J^0 \partial J^{n^2+2n}:$. We may disregard $:J^0 \partial J^{n^2+2n}:$ since it lies in $C_2 (\mathcal{M}_{-n})$, and hence will not contribute to the symbol of $\pi_{Zh}(E)$ in $A(\mathcal{M}_{-n})$, which has weight $(n+1)^2 + 2$. Finally, the terms $\partial^2 J^{n^2+2n}$ and $\partial J^{n^2+2n+1}$ which can appear in $E^2$, may be disregarded as well, since they lie in $C_2(\mathcal{M}_{-n})$.

Using the relation $D_{(0,\dots,n),(0,\dots,n-1,n+1)}$, we can eliminate $J^{n^2+2n+1}$ and its derivatives from $E$, since $R_{(0,\dots,n),(0,\dots,n-1,n+1)}$ is nonzero. Next, we can eliminate $J^{n^2+2n}$ and its derivatives using $D_0 = D_{(0,\dots,n),(0,\dots,n)}$, since $R_0\neq 0$. Let $E'$ be the element of $\mathcal{I}_{-n}$ obtained from $E$ in this way. Since $E'$ only depends on $J^0,\dots,J^{n^2+2n-1}$, $\pi_{Zh}(E') \in A(\mathcal{M}_{-n})$ will lie in $\mathbb{C}[a^0,\dots,a^{n^2+2n-1}]$. We will see that $\pi_{Zh}(E')$ is nonzero, and hence gives rise to a nontrivial relation among the generators $a^0,\dots,a^{n^2+2n-1}$ of $A( \mathcal{W}_{1+\infty,-n})$. It is enough to show that the symbol of $\pi_{Zh}(E')$ in $gr(A(\mathcal{M}_{-n}))$ has weight $(n+1)^2+2$, and is nonzero.

We introduce the {\it degree-lexicographic monomial ordering} on $\mathbb{C}[a^0,a^1,a^2,\dots]$, where the variables are ordered by $a^l < a^{l+1}$, for $l\geq 0$. Given a polynomial $P\in \mathbb{C}[a^0,a^1,a^2,\dots]$, $Symb(P)$ will denote the component of maximal weight, and $LT(P)$ will denote the leading term of $Symb(P)$. Similarly, given a vertex operator $F\in \mathcal{M}_{-n}$, let $LT(F)$ denote $LT(\pi_{Zh}(F))$, i.e., the leading term of the symbol of $\pi_{Zh}(F)$. An easy calculation shows that $\Omega_{k,l} \cong (-1)^{l} J^{k+l}$ modulo total derivatives, so $LT (\Omega_{k,l}) = (-1)^l a^{k+l}$. 
It follows that $LT(D_0) = \pm \prod_{k=0}^n a^{2k}$. Similarly, we have $$LT(D_{(0,\dots,n),(0,\dots,n-1,n+2)})= \pm a^{2n+2} \prod_{k=0}^{n-1} a^{2k},$$ $$LT(D_{(0,\dots,n),(0,\dots,n-2,n,n+1)}) = \pm a^{2n-1}a^{2n+1} \prod_{k=0}^{n-2} a^{2k}.$$

Now we return to our element $E$ given by (\ref{defofe}). Suppose first that $E^4$ contains the term $:J^0 J^{n^2+2n+1}:$ with nonzero coefficient. As above, using the relation $D_{(0,\dots,n),(0,\dots,n-1,n+1)}$, we eliminate the variable $J^{n^2+2n+1}$, and using $D_0$, we eliminate $J^{n^2+2n}$, to obtain $E'\in \mathcal{I}_{-n}$. Since $E^4$ contains $:J^0 J^{n^2+2n+1}:$ with nonzero coefficient, and $D^4_{(0,\dots,n),(0,\dots,n-1,n+1)}$ contains $:J^0 J^{n^2+2n}:$ with nonzero coefficient, it follows that $E'$ contains the homogeneous, normally ordered polynomial $:J^0 J^0 D_0^{2n+2}:$ in the variables $J^0,\dots,J^{n^2+2n-1}$ and their derivatives. (Recall that $D^{2n+2}_0$ is the term of maximal degree appearing in the decomposition $D_0 = \sum_{k=1}^{n+1} D^{2k}_0$). It follows that $LT(E') = \pm (a^0)^2 \prod_{k=0}^n a^{2k}$, and hence is nonzero.

Next, suppose that $E^4$ does not contain $:J^0 J^{n^2+2n+1}:$, but that it does contain the term $:J^1 J^{n^2+2n}:$ with nonzero coefficient. A similar argument shows that $LT(E') = \pm a^1 \prod_{k=0}^{n} a^{2k}$, and hence is nonzero. Finally, suppose that $E^4$ contains neither $:J^0 J^{n^2+2n+1}:$ nor $:J^1 J^{n^2+2n}:$ with nonzero coefficient. Then the symbol of $\pi_{Zh}(E')$ has degree at most $2n+2$, and coincides with the symbol of $\pi_{Zh}(E)$. It follows that $$LT(E') =  LT(D_{(0,\dots,n),(0,\dots,n-1,n+2)}) = \pm  a^{2n+2} \prod_{k=0}^{n-1} a^{2k}.$$ In particular, $LT(E')$ is nonzero, as desired. \end{proof}

It is an interesting problem to calculate the dimension of the variety $V(I_{-n})$ and determine whether it is irreducible or not. We hope to return to these questions in future work.

\end{document}